\documentclass[11pt]{amsart} 
\calclayout

\usepackage{setspace} 
\usepackage{graphicx}
\usepackage{microtype}

\usepackage{pgfplots}
\pgfplotsset{compat=1.15}
\usepackage{mathrsfs}

\usepackage{mathrsfs}

\usetikzlibrary{arrows}

\usepackage[english]{babel}
\usepackage[maxbibnames=5,minbibnames=4,backend=bibtex]{biblatex}
\usepackage{csquotes}
\usepackage{amssymb}
\usepackage{amsmath}
\usepackage{amsthm}

\newtheorem{thm}{Theorem}[section]
\newtheorem{lem}[thm]{Lemma}
\newtheorem{prop}[thm]{Proposition}
\newtheorem{coro}[thm]{Corollary}

\theoremstyle{remark}
\newtheorem{rema}[thm]{Remark}
\newtheorem{exa}[thm]{Example}
\newtheorem{defi}[thm]{Definition}

\DeclareFontFamily{U} {MnSymbolA}{}
\DeclareFontShape{U}{MnSymbolA}{m}{n}{
  <-6> MnSymbolA5
  <6-7> MnSymbolA6
  <7-8> MnSymbolA7
  <8-9> MnSymbolA8
  <9-10> MnSymbolA9
  <10-12> MnSymbolA10
  <12-> MnSymbolA12}{}
\DeclareFontShape{U}{MnSymbolA}{b}{n}{
  <-6> MnSymbolA-Bold5
  <6-7> MnSymbolA-Bold6
  <7-8> MnSymbolA-Bold7
  <8-9> MnSymbolA-Bold8
  <9-10> MnSymbolA-Bold9
  <10-12> MnSymbolA-Bold10
  <12-> MnSymbolA-Bold12}{}

\DeclareSymbolFont{MnSyA} {U} {MnSymbolA}{m}{n}
\DeclareMathSymbol{\lcirclearrowright}{\mathrel}{MnSyA}{252}
\DeclareMathSymbol{\rcirclearrowleft}{\mathrel}{MnSyA}{250}

\newcommand{\rank} 
{\operatornamewithlimits{rank}}

\usepackage{stmaryrd}
\usepackage{slashed}
\usepackage{tikz}
\usetikzlibrary{positioning}
\usetikzlibrary{calc,intersections}
\usepackage{todonotes}
\usepackage[all,cmtip]{xy}
    
\title{Transitive Courant algebroids and double symplectic groupoids}
\author{Daniel \'Alvarez}\address{Departament of Mathematics, University of Toronto, Toronto, Ontario}
\email{dalv@math.toronto.edu}

\date{} 
\begin{document}

\begin{abstract} In this work we extend the Lu-Weinstein construction of double symplectic groupoids to any Lie bialgebroid such that its associated Courant algebroid is transitive under natural integrability assumptions. We illustrate this result by showing how it generalises many of the examples of double symplectic groupoids that have appeared in the literature. As preliminary steps for this construction, we give a classification of exact twisted Courant algebroids over Lie groupoids (CA-groupoids for short) and we show the existence of a foliation by twisted Courant algebroids on the base of a twisted CA-groupoid. 
\end{abstract}\maketitle
\tableofcontents 

\section{Introduction} Symplectic groupoids are fundamental tools in Poisson geometry, being the global counterparts of Poisson structures \cite{weisygr,craruimar}. Not every Poisson manifold can be associated with a symplectic groupoid, only those that satisfy the necessary and sufficient conditions of \cite{craruipoi} admit such a corresponding global object. Interestingly, some Poisson structures that are compatible with a group(oid) structure such as those defined on Poisson-Lie groups give rise to symplectic groupoids that are equipped with an additional Lie groupoid structure which makes them into {\em double symplectic groupoids} \cite{luwei2}. Some of the most interesting applications of double symplectic groupoids are the following. In the case of a Poisson-Lie group, the presence of an additional groupoid multiplication on an associated symplectic groupoid is responsible for the existence of a monoidal structure on the category of its hamiltonian spaces, such  a monoidal structure is even braided monoidal in some situations \cite{weixuybe}; we expect that similar properties should hold for more general double symplectic groupoids by using one of the groupoid multiplications to define the monoidal structure. Another remarkable aspect of double symplectic groupoids is their relevance to the AKSZ formalism \cite{aksz}. The infinitesimal counterpart of a double symplectic groupoid is a Lie bialgebroid \cite{macxu} and the doubles of Lie bialgebroids are Courant algebroids \cite{liuweixu}. This means that double symplectic groupoids induce integrations of Courant algebroids \cite{rajtan,cuezhu} and so they can be regarded as classical phase spaces associated with the Courant sigma model \cite{royaksz}. It was pointed out in \cite{stelagpd} (after \cite{sevsomtit}) that this viewpoint could lead to a full integration of Lie bialgebroids; conversely, we expect that the explicit double symplectic groupoids we construct can be used to better understand the Courant sigma model in special circumstances, see \cite{cal:lag}. 

The construction of double symplectic groupoids has proceeded so far on a case by case basis without a general principle behind. In this work we describe a way of constructing double symplectic groupoids that can be applied to most cases of interest, motivated by the ideas in \cite{luwei2,sevmorqua,sevlet}. Instead of attempting to integrate a Lie bialgebroid $(A,B)$ by stages, first to a Poisson groupoid \cite{macxu2} and then to a double symplectic groupoid \cite{stelagpd}, we propose to look at both $A$ and $B$ as {\em Dirac structures} and integrate them according to their interpretation as {\em infinitesimal 2-shifted lagrangians} \cite{pymsaf}; these 2-shifted lagrangians can then be used to produce a double symplectic groupoid following the outline of \cite{luwei2}. Under mild assumptions, our approach allows us to integrate a Lie bialgebroid whose associated Courant algebroid is transitive to a double symplectic groupoid, thus giving a unified explanation for most of the known examples such as those coming from dynamical Poisson groupoids \cite{dynpoi} or moduli spaces of flat bundles \cite{luwei2,sevmorqua}. We can also use this construction to produce new families of double symplectic groupoids in two interesting contexts: (1) in generalised K\"ahler geometry, where the relevant Lie bialgebroids are holomorphic and induce exact Courant algebroids \cite{guagk} and (2) in Lie theory, where we can give integrations in the algebraic category of the configuration Poisson groupoids of flags \cite{conpoigro}, these applications will be given in upcoming works. 

Let us describe the contents of this paper in more detail. We start with a classification of exact Courant algebroid groupoids (CA-groupoids) in terms of certain degree 4 cohomology classes in the Bott-Shulman-Stasheff complex of a Lie groupoid, see Theorem \ref{thm:cla}; this is analogous to the classification of exact Courant algebroids over manifolds of \cite{sevlet}. Next we study transitive Courant algebroids from a Lie groupoid perspective shedding new light on some of their well known features, see Theorem \ref{thm:hetca} which is a global counterpart of \cite[Corollary 5.10]{pymsaf}. A transitive Courant algebroid ${E}$ over a manifold $M$ determines a transitive Lie algebroid $E/T^*M$ which corresponds, if it is integrable, to the gauge groupoid $\mathcal{G}=(P\times P)/G \rightrightarrows M$ of a principal $G$-bundle $P\rightarrow M$. The Courant algebroid $E$ acts naturally on $\mathcal{G} \rightrightarrows M$ inducing an exact CA-groupoid $\mathtt{t}^*E\oplus \mathtt{s}^*E \rightrightarrows E$ over it. According to our classification of exact CA-groupoids, a splitting of $\mathtt{t}^*E\oplus \mathtt{s}^*E \rightrightarrows E$ determines a simplicially closed 2-form $\omega_2$ defined on the space of composable pairs in $\mathcal{G}$; $\omega_2$ is {\em closed up to homotopy} in the sense that it is de Rham closed up to the simplicial differential of a closed 3-form $\omega_1$ on $\mathcal{G}$, moreover $\omega_2$ is {\em nondegenerate} in the sense of shifted symplectic geometry \cite{ptvv,getdia}, see Definition \ref{def:2shisym}. This means that $\omega_2+\omega_1$ is a 2-shifted symplectic structure. Following \cite{pymsaf,ptvv} and based on the previous observation, a Dirac structure $A\subset E$ is integrable by a Lie groupoid $\mathcal{L}(A) \rightrightarrows M$ equipped with a Lie groupoid morphism $\Phi_A:\mathcal{L}(A) \rightarrow \mathcal{G}$ and $\beta_A\in \Omega^2(\mathcal{L}(A))$ that serves as a primitive for the pullback of $\omega_2+\omega_1$ along $\Phi_A$ in the Bott-Shulman-Stasheff complex of $\mathcal{L}(A) \rightrightarrows M$. We shall call $(\mathcal{L}(A),\Phi_A,\beta_A)$ a {\em 2-shifted lagrangian groupoid}. If $B\subset E$ is another Dirac structure which is complementary to $A$ and there is a corresponding 2-shifted lagrangian groupoid $(\mathcal{L}(B),\Phi_B,\beta_B)$ that integrates $B$, we can construct a double symplectic groupoid with side groupoids $\mathcal{L}(A)$ and $\mathcal{L}(B)$ that integrates the Lie bialgebroid $(A,B)$, see Theorem \ref{thm:dousym}. 

This work is complementary to \cite{2lagpoi} where we describe how to perform the integration of Dirac structures in transitive Courant algebroids to 2-shifted lagrangian groupoids as well as many explicit examples, here we shall describe how to produce a double symplectic groupoid as above assuming that the corresponding 2-shifted lagrangian groupoids are given. Along the way, we make some general observations that can be of independent interest: we show that there is a foliation by Courant algebroids on the base of a CA-groupoid in Theorem \ref{thm:folca}, we introduce the concept of double quasi-symplectic groupoid and we show that it is the global counterpart of a multiplicative Dirac structure (Proposition \ref{pro:difdir}), finally, we provide in Theorem \ref{thm:2sym} an extension to double quasi-symplectic groupoids of the correspondence established in \cite{rajtan} between double symplectic groupoids and certain local 2-shifted symplectic Lie 2-groupoids. We use shifted symplectic geometry \cite{ptvv,getdia,pymsaf} as conceptual guidance but we only rely on Lie theory and Poisson geometry.

\section{Classification of exact twisted CA-groupoids}
\subsection{Preliminaries} Let $H$ be a closed 4-form on a manifold $M$. An $H$-twisted Courant algebroid \cite{hanstr} is a vector bundle $E \rightarrow M$ equipped with a non-degenerate symmetric bilinear pairing $\langle \,,\, \rangle $, a vector bundle map $\mathtt{a}:E \rightarrow TM$ and a bilinear bracket $[\,,\,]$ on $\Gamma (E)$ such that
\begin{align*} &\mathtt{a}^*d \langle e_1,e_2 \rangle =[e_1,e_2]+[e_2,e_1] \\
  &[e_1,fe_2]=f[e_1,e_2]+\mathcal{L}_{\mathtt{a}(e_1) }(f )e_2 \\
& \mathcal{L}_{\mathtt{a}(e_1) }\langle e_2,e_3 \rangle =\langle [e_1,e_2],e_3 \rangle + \langle e_2,[e_1,e_3] \rangle \\ 
 &[e_1,[e_2,e_3]]=[[e_1,e_2],e_3]+[e_2,[e_1,e_3]]+\mathtt{a}^*(i_{\mathtt{a}(e_1) }i_{\mathtt{a}(e_2) }i_{\mathtt{a}(e_3) }H)
\end{align*} 
for all $e_i\in \Gamma (E)$ and all $f\in C^\infty(M)$. From these equations it follows that $\mathtt{a}\circ \mathtt{a}^*=0 $ and that $\mathtt{a}$ preserves brackets, among other identities \cite{hanstr}. We call $H$ the twisting form of $E$. Take any 3-form $C$ on $M$, the basic example for us is $\mathbb{T}_{ C}M=TM \oplus T^*M$ equipped with the canonical pairing, denoted by $\langle \,,\, \rangle $, and the projection to $TM$ as the anchor; the bracket is the Courant-Dorfman bracket
\[ [X\oplus \alpha ,Y \oplus \beta ]=[X,Y]\oplus \mathcal{L}_X \beta -i_Y d \alpha -i_Xi_Y C; \]
we have that $\mathbb{T}_{C}M$ is a $dC$-twisted Courant algebroid which is exact, meaning that the anchor sequence defined by $\mathtt{a}$ and $\mathtt{a}^*$ is exact. Just as in the untwisted case, every exact twisted CA-groupoid is of this form. If $\mathbb{T}_C M$ is an exact $dC$-twisted Courant algebroid and $\omega \in \Omega^2(M)$, the $B$-field (or gauge) transformation induced by $\omega $ is the map $e^{\omega }(X\oplus \alpha )=X\oplus \alpha +i_{X}\omega $ and it constitutes an isomorphism $\mathbb{T}_C M \cong \mathbb{T}_{C-d \omega }M$.   

If $E \rightarrow M$ is an $H$-twisted Courant algebroid and we take a submanifold $i:N \hookrightarrow M$ such that $i^*H$ is exact, a twisted Dirac strucure with support $N$ is a vector subbundle $L $ over $N$ that satisfies $L=L^\perp$ and such that the space of sections of $E$ that restrict to sections of $L$ is closed under the restricted bracket \cite{pymsaf}. An {\em exact Dirac structure} is a twisted Dirac structure for which the restricted anchor map is onto the tangent bundle of the support manifold; these objects are also known as isotropic involutive splittings \cite{baigua}.

A Lie groupoid is a groupoid object in the category of smooth manifolds, denoted by $G_1 \rightrightarrows G_0$, such that its source map is a submersion \cite{macgen}. The structure maps of a Lie groupoid are its source, target, unit, multiplication and inversion maps, denoted respectively $\mathtt{s}$, $\mathtt{t}$, $\mathtt{u}$, $\mathtt{m}$, $\mathtt{i}$; in order to to avoid ambiguities we put a subindex next to the structure maps of a Lie groupoid when necessary. We will also use the abbreviations $\mathtt{m}(a,b)=ab$ and $\mathtt{i}(a)=a^{-1}$ to simplify the notation. We denote by $A_{G_1}=\ker T\mathtt{s}|_{G_0}$ the Lie algebroid of $G_1\rightrightarrows G_0$ \cite{macgen}. 

A VB-groupoid is a Lie groupoid which is also a groupoid in the category of vector bundles, so that all its structure maps are vector bundle morphisms \cite{macgen}. 

\begin{defi} Let $G_1 \rightrightarrows G_0$ be a Lie groupoid. A {\em twisted CA-groupoid over $G_1 \rightrightarrows G_0$} is a VB-groupoid $\mathbb{E}_1 \rightrightarrows \mathbb{E}_0$ over $G_1 \rightrightarrows G_0$ such that:
\begin{enumerate} \item $\mathbb{E}_1$ is a $(\mathtt{s}^*\omega_0 -\mathtt{t}^* \omega_0)$-twisted Courant algebroid for some $\omega_0\in \Omega^4(G_0)$ closed; 
\item the graph of the multiplication in $\mathbb{E}_1$ is a Dirac structure 
\[ \text{graph}(\mathtt{m}_{\mathbb{E}_1 } )=\{(\mathtt{m}_{\mathbb{E}_1 }(u,v),u,v)|\mathtt{s}_{\mathbb{E}_1 }(u)=\mathtt{t}_{\mathbb{E}_1 }(v)\} \hookrightarrow  \overline{\mathbb{E}_1 } \times \mathbb{E}_1 \times \mathbb{E}_1  \]
with support the graph of the multiplication in $G_1 \rightrightarrows G_0$; the notation $\overline{\mathbb{E}_1 }$ means that we are changing the sign of the metric of this factor. \end{enumerate} 
In this situation, we say that $\mathbb{E}_1 \rightrightarrows \mathbb{E}_0$ is an $\omega_0$-twisted CA-groupoid and $\omega_0$ is called its {\em multiplicative twisting form}. Suppose that $(H_1 \rightrightarrows H_0) \hookrightarrow (G_1 \rightrightarrows G_0)$ is a Lie subgroupoid such that the pullback of $\omega_0$ is exact on $H_0$; a {\em twisted multiplicative Dirac structure $L_1 \rightrightarrows L_0$} inside $\mathbb{E}_1 \rightrightarrows \mathbb{E}_0$ with support $H_1 \rightrightarrows H_0$ is a VB-subgroupoid such that the pullback of $\omega_0$ to $H_0$ is exact and $L_1\hookrightarrow \mathbb{E}_1$ is a twisted Dirac structure supported on $H_1$. \end{defi} 
\begin{rema} Notice that for an exact twisted CA-groupoid the twisting form on the space of arrows of the underlying Lie groupoid is automatically exact. \end{rema} 
We introduce this definition in order to realize geometrically the cohomology classes of certain degree four elements in the Bott-Shulman-Stasheff complex of a Lie groupoid, the usual concept of CA-groupoid is not sufficiently general as we will see later on.
  
If $G_1 \rightrightarrows G_0$ is a Lie groupoid, we denote by $G_n$ the manifold of all the $n$-tuples of composable elements in $G_1$, the projections, multiplication and unit inclusions make this collection of spaces into a simplicial manifold $(G_\bullet,d_i,s_i) $ called the {\em nerve of $G_1 \rightrightarrows G_0$} which we also denote just by $G_\bullet$. Explicitly:
\begin{align*} &G_n=\{(x_1,\dots,x_n)|\mathtt{s}(x_i)=\mathtt{t}(x_{i+1}), \, i=1\dots n-1  \}, \\ 
&d_1= \mathtt{t},\, d_0=\mathtt{s}:G_1 \rightarrow G_0\\
& d_i(x_1,\dots,x_n)=\begin{cases} (x_2,\dots,x_n) \quad \text{if $i=0$}, \\ 
(x_1,\dots,x_{n-1})\quad  \text{if $i=n$}, \\ 
(x_1,\dots,x_{i-1},x_ix_{i+1},x_{i+2},\dots, x_n), \quad  \text{if $ 1\leq i \leq n-1$}; \end{cases}  \\   
&s_i(x_1,\dots,x_n)=(x_1,\dots,x_{i},\mathtt{u}({\mathtt{t}(x_{i+1}) }), x_{i+1},x_{i+2},\dots, x_n)\quad  0\leq i \leq n-1, \\ 
&s_n(x_1,\dots,x_n)=(x_1,\dots,x_n,\mathtt{u}(\mathtt{s}(x_n)))
&\forall (x_1,\dots,x_n)\in G_n. \end{align*} 

The spaces of differential forms $\Omega^a(G_b) $ form a double complex $ \Omega^\bullet(G_\bullet)$, called the Bott-Shulman-Stasheff complex, in which we have the de Rham differential $d$ and the simplicial differential $\delta=\sum_i (-1)^id_i^*$. So the associated total complex $\text{Tot}_n(\Omega^\bullet(G_\bullet))=\bigoplus_{a+b=n}\Omega^a(G_b)$ is equipped with a differential $D$ which is $D=(-1)^bd+\delta$ restricted to each $\Omega^a(G_b)$. An $m$-shifted $n$-form $\omega $ on $G_1 \rightrightarrows G_0$ is a collection of forms $\omega_i\in \Omega^{m+n-i}(G_i) $ with $i=0\dots m$ \cite{getdia}. A form $\omega $ is {\em closed} if $D(\sum_i \omega_i)=0$ and {\em normalised} if $s^*_i \omega_j=0$ for all $i$ and $j$. Notice that these definitions make sense for every simplicial manifold. The space of normalised forms $\widehat{\Omega }^\bullet(G_\bullet)$ constitutes a subcomplex of $(\text{Tot}_n(\Omega^\bullet(G_\bullet)))_n$. The complex of $n$-forms is 

\[ \mathcal{A}^n ( G_\bullet)_m=\bigoplus_{a\geq n,\, b\geq0, \, a+b=n+m} \widehat{\Omega}^a(G_b),  \quad m\geq 0  \] 
with $D$ as the differential. We will only be interested in the complex $\mathcal{A}^2 (G_\bullet)$, the cohomology  spaces of $\mathcal{A}^2 (G_\bullet)$ can be interpreted as the 2-forms on the differentiable stack represented by $G_1 \rightrightarrows G_0$ with the shift of a 2-form given by its degree as a cohomology class \cite{getdia}.

\subsection{The classification} Let $\mathbb{E}_1 \rightrightarrows \mathbb{E}_0 $ be an exact CA-groupoid over a Lie groupoid $G_1 \rightrightarrows G_0$. The choice of a splitting of the anchor sequence 
\[ \xymatrix{0 \ar[r] & T^*G_1 \ar[r] & \mathbb{E}_1 \ar[r] & TG_1 \ar[r] & 0 } \]
allows us to identify $\mathbb{E}$ with the $\omega_1$-twisted Courant algebroid $\mathbb{T}_{\omega_1}G_1$ for some $\omega_1\in \Omega^3(G_1)$. The canonical (untwisted) Courant algebroid $\mathbb{T}G_1$ is a CA-groupoid \cite{muldirphd}. We can see that the graph of the multiplication $\mathtt{m}_{\mathbb{E}_1 } $ on $\mathbb{E}_1$ corresponds under this identification to an exact Dirac structure inside $ \overline{\mathbb{T}G_1} \times  \mathbb{T}G_1 \times \mathbb{T}G_1$ which is the result of applying a B-field transformation to the canonical multiplication of $\mathbb{T}G_1$. In other words, it is the graph of some $\omega_2\in \Omega^2(G_2)$ over the graph of the multiplication on $G_1 \rightrightarrows G_0$ and this graph completely determines the groupoid structure of $\mathbb{E}_1$. 

In this first proposition we will see that a normalised 2-shifted 2-form determines a deformed groupoid structure on $\mathbb{T}G_1$ and the effect of choosing a different splitting of this Courant algebroid corresponds to the action of an exact element $D \eta_1$ for $\eta_1\in \widehat{\Omega }^2(G_1)$ as a B-field transformation. On Theorem \ref{thm:cla} we will see that two closed 2-shifted 2-forms that differ by $D\eta_0$ for $\eta_0\in \Omega^3(G_0)$ induce twisted CA-groupoids which are Morita equivalent in a suitable sense and this completes the classification.

\begin{prop}\label{thm:excagpd} The choice of a splitting of the anchor sequence establishes a bijection between isomorphism classes of exact twisted CA-groupoids over $G_1 \rightrightarrows G_0$ up to a base fixing automorphism and the space
\[   (\ker D  \cap \mathcal{A}^2 (G_\bullet)_2)/D \widehat{\Omega}^2 (G_1). \]
\end{prop} 
    
\begin{proof} {\em Part 1: normalised closed 2-shifted 2-forms determine twisted CA-groupoid structures}. Take $\omega =\omega_0+\omega_1+\omega_2\in \mathcal{A}^2 (G_\bullet)_2$ such that $D \omega =0$. Let us spell the equations that these forms satisfy: 
\begin{align*} d\omega_0=0, \quad d \omega_1=\delta \omega_0, \quad d \omega_2 =-\delta \omega_1, \quad \delta \omega_2=0.  \end{align*}   
Denote by $A=A_{G_1}$ the Lie algebroid of $G_1 \rightrightarrows G_0$. We can define a CA-groupoid structure on the twisted Courant algebroid $\mathbb{E}_1:= \mathbb{T}_{\omega_1 }G_1$ over $\mathbb{E}_0:=TG_0\oplus A^*$ as follows. Let $p_T:\mathbb{E}_1 \rightarrow TG_1$ and $p_{T^*}:\mathbb{E}_1 \rightarrow T^*G_1$ be the respective projections. In order to specify how $\omega_2$ deforms the canonical groupoid structure on $\mathbb{T}G_1$, since $d \omega_2 =-\delta \omega_1$, we can consider the following Dirac structure which is supported on $\text{graph} (\mathtt{m}_{G_1})$:
\begin{align} L_{\omega }=\left\{ (X\oplus \alpha ,Y\oplus \beta ,Z\oplus \gamma) \in \overline{\mathbb{T}G_1}\times\mathbb{T}G_1\times \mathbb{T}G_1\left| \begin{aligned} 
&\mathtt{m}^*_{G_1} \alpha =\text{pr}_1^*\beta + \text{pr}_2^*\gamma -i_{(Y,Z)}\omega_2 \\
& X=T \mathtt{m}(Y,Z), \end{aligned}   \right.\right\}; \label{eq:twimul}
    \end{align} 
note that this Dirac structure can be viewed as an extension of the graph of the canonical multiplication on $(T^*G_1\rightrightarrows A^*)\hookrightarrow (\mathbb{E}_1\rightrightarrows \mathbb{E}_0)$ and for $\omega_2=0$ it is the graph of the canonical multiplication on $\mathbb{T}G_1$. This Dirac structure determines completely a new groupoid structure on $\mathbb{E}_1\rightrightarrows \mathbb{E}_0$; we shall denote the corresponding structure maps as $\mathtt{s}_{\omega },\mathtt{t}_{\omega },\mathtt{u}_{\omega },\mathtt{m}_{\omega },\mathtt{i}_{\omega }$, see Proposition \ref{pro:twimul} for their explicit description.

{\em Claim: the twisted CA-groupoids associated to representatives of the same class in $(\ker D  \cap \mathcal{A}^2 (G_\bullet)_2)/D \widehat{\Omega}^2 (G_1)$ are related by a $B$-field transformation}. If $\omega $ and $\omega'$ are $D$-closed forms such that $\omega -\omega'=D \theta $, where $\theta \in \Omega^2(G_1)$ is normalised as well, we have that the $B$-field transformation $\tau$: $X \oplus \alpha \mapsto X\oplus \alpha +i_X \theta $ is a twisted Courant algebroid isomorphism $\mathbb{T}_{\omega'_1}G_1 \rightarrow \mathbb{T}_{\omega'_1-d\theta }G_1  $ that also satisfies $ \mathtt{t}_{\omega }\circ\tau=\mathtt{t}_{\omega'}   $, $ \mathtt{s}_{\omega }\circ\tau=\mathtt{s}_{\omega'}   $ and hence $\tau$ takes $L_{\omega'}$ to $L_{\omega'+\delta \theta }=L_{\omega}$ which implies that it is a groupoid isomorphism.
   
{\em Part 2: splitting an exact twisted CA-groupoid determines a normalised closed 2-shifted 2-form on the base groupoid}. Given an exact $\omega_0$-twisted CA-groupoid $\mathbb{E}_1 \rightrightarrows \mathbb{E}_0$ over $G_1 \rightrightarrows G_0$, splitting the anchor sequence allows us to identify $\mathbb{E}_1$ with $\mathbb{T}_{\omega_1}G_1$ for some $\omega_1\in \Omega^3(G_1)$ such that $d\omega_1=\delta \omega_0$ and so we can also identify $\mathbb{E}_0 $ with $A^* \oplus TG_0$; since $\mathbb{E}_0\hookrightarrow \mathbb{E}_1$ is an exact Dirac structure supported on $G_0$, we can take a splitting such that $\mathtt{u}^*_{G_1}\omega_1=0$. By computing the rank of the graph $L \hookrightarrow \overline{ \mathbb{E} }_1\times \mathbb{E}_1 \times \mathbb{E}_1  $ of $\mathtt{m}_{\mathbb{E}_1 } $ and noticing that the kernel of the anchor map restricted to $L$ is the graph of the multiplication on $T^*G_1$, we get that $L$ is a Dirac structure over $\text{graph}(\mathtt{m}_{G_1}) $ which covers $T\text{graph}(\mathtt{m}_{G_1}) $ surjectively and hence it has to be of the form $L_{\omega }$ as above for some $\omega_2\in \Omega^2(G_2)$ that has to satisfy $d \omega_2= -\delta \omega_1$. As we show in the proof of Proposition \ref{pro:twimul}, the fact that $\mathtt{m}_{\mathbb{E}_1 } $ is associative implies that $\delta \omega_2=0$ and the normalisation follows from the unit axiom.  \end{proof}
In what follows we shall denote the CA-groupoid corresponding to a closed 2-shifted 2-form $\omega=\omega_0+\omega_1 + \omega_2$ on $G_1 \rightrightarrows G_0$ as $(\mathbb{T}G_1 \rightrightarrows TG_0\oplus A_{G_1}^*, \omega ) $, where $A_{G_1}=\ker T\mathtt{s}|_{G_0} $ is the Lie algebroid of $G_1 \rightrightarrows G_0$.
\begin{coro}\label{cor:exdir} Let $G_1 \rightrightarrows G_0$ be a Lie groupoid and let $i:(H_1 \rightrightarrows H_0) \hookrightarrow (G_1 \rightrightarrows G_0)$ be a Lie subgroupoid. If $\omega =\omega_1+ \omega_2$ is a closed 2-shifted 2-form on $G_1 \rightrightarrows G_0$, then the exact multiplicative Dirac structures inside the exact CA-groupoid $(\mathbb{T}G_1 \rightrightarrows TG_0\oplus A_{G_1}^*,\omega ) $ with support $H_1 \rightrightarrows H_0$ are classified by normalised primitives $\beta_1\in \Omega^2(H_1)$ for the pullback $i^*\omega $ in the Bott-Shulman-Stasheff complex of $H_1 \rightrightarrows H_0$. \end{coro} 
\begin{proof} An exact Dirac structure in $(\mathbb{T}G_1 \rightrightarrows TG_0\oplus A_{G_1}^*,\omega ) $ with support $H_1 \rightrightarrows H_0$ has to be by definition the graph of a 2-form $\beta_1\in \Omega^2(H_1) $ such that $d \beta_1 =-i^*\omega_1$. But such a Dirac structure is multiplicative if and only if
\[ \langle i_{T \mathtt{m}(X,Y) } \beta_1, T \mathtt{m}(U,V) \rangle =\langle i_X \beta_1, U \rangle + \langle i_Y \beta_1,V \rangle - \langle i_{(X,Y)} \omega_2 ,(U,V)  \rangle \] 
for all $(U,V),(X,Y)\in TH_2$ and this just means that $\delta \beta_1= i^*\omega_2$. From the same equation applied to $U=V$ and $X=Y\in TG_0$ we get that $\beta_1$ is normalised. \end{proof} 
\begin{exa}\label{exa:bg} Let $K$ be a Lie group such that its Lie algebra $\mathfrak{k}$ is quadratic: i.e. it is endowed with an Ad-invariant nondegenerate symmetric bilinear form $\langle \,,\, \rangle $ and let $k^{-1}dk $ and $(dk)k^{-1}$ denote, respectively, the left and right invariant Maurer-Cartan forms on $K$ at $k$. Then $\Omega=\Omega_1+\Omega_2$ given by
\begin{align} \Omega_1|_k=\frac{1}{6}\langle k^{-1}dk\wedge [k^{-1}dk\wedge k^{-1}dk]\rangle, \quad \Omega_2|_{(x,y)}=-\langle x^{-1}dx \wedge (dy)y^{-1} \rangle\quad \forall k,x,y\in K \label{eq:bg} \end{align}   
is a closed 2-shifted 2-form $\Omega $ on $K_\bullet$. The CA-groupoid it induces on $K \rightrightarrows \{1\}$ is the one that appears in \cite{purspi} as an action Courant algebroid in which the Cartan-Dirac structure lies. See also \cite{liedir} for other versions of this example. 

There is a classical variant of the Chern-Weil homomorphism due to Shulman \cite{shuphd} that associates closed shifted $n$-forms on $K_\bullet$ to invariant polynomials on $\mathfrak{k}^*$, one can check that Shulman's map takes $P=\langle\,,\, \rangle$ to a multiple of $\Omega$. Moreover, in the case of $K=SU(n)$ and the pairing given by the trace, Shulman's map applied to $P$ gives us a representative of the classical universal first Pontryagin class \cite{getdia}. For this reason, $\Omega$ can be called a representative of the generalised universal first Pontryagin class associated with $\langle\,,\, \rangle$. See also \cite[\S 4.2]{walmulger}. \end{exa} 

Now we shall introduce the notion of Morita equivalence for twisted CA-groupoids. Recall that \emph{a left (right) Lie groupoid action} of a Lie groupoid $G_1 \rightrightarrows G_0$ on a map $J:S \rightarrow G_0$ is a smooth map $a:(G_1)_\mathtt{s} \times_J S \rightarrow S$ (respectively, $a:S_J \times_\mathtt{t} G_1 \rightarrow S$) which satisfies (1) $a(\mathtt{m} (g,h),x)=a(g,a(h,x))$ (respectively, $a(x,\mathtt{m} (g,h))=a(a(x,g),h)$) for all $g,h\in G_1$ and for all $x\in S$ for which $a$ and $\mathtt{m}$ are defined and (2) $a(\mathtt{u}(J(x)),x)= x$ for all $x\in G_0$. For the sake of brevity, we denote the fiber product $(G_1)_\mathtt{s} \times_J S $ by $G_1 \times_{G_0} S$ and $a(g,x)=g\cdot x$. In this situation, $G_1 \times_{G_0} S$ becomes a Lie groupoid over $S$ with the projection being the source map, $a$ being the target map and the multiplication is defined by $(g,a(h,p))(h,p)=(gh,p)$. The resulting Lie groupoid is denoted by $G_1 \ltimes_{G_0} S\rightrightarrows S$ (respectively, $S \rtimes_{G_0} G_1 \rightrightarrows S$) and it is called an {\em action groupoid}.
\begin{defi} Let $\mathbb{E}_1 \rightrightarrows \mathbb{E}_0$ and $\mathbb{F}_1 \rightrightarrows \mathbb{F}_0$ be two twisted CA-groupoids over $G_1 \rightrightarrows G_0$ and $H_1 \rightrightarrows H_0$ respectively with $\omega_0 \in \Omega^4(G_0)$ and $\omega_0'\in \Omega^4(H_0)$ as the corresponding multiplicative twisting forms. A {\em Courant Morita equivalence} between $\mathbb{E}_1 \rightrightarrows \mathbb{E}_0$ and $\mathbb{F}_1 \rightrightarrows \mathbb{F}_0$ consists of the following data. Take a twisted Courant algebroid $\mathbb{A}$ over $S$ with exact twisting form $dH$ equipped with two vector bundle maps $p:\mathbb{A} \rightarrow \mathbb{E}_0 $ and $q:\mathbb{A} \rightarrow \mathbb{F}_0 $ over $p_0:S \rightarrow G_0$ and $q_0:S \rightarrow H_0$ together with two Lie groupoid actions over those maps
\begin{align*} \mathbb{E}_1 \times_{\mathbb{E}_0 } \mathbb{A} \rightarrow \mathbb{A} \\
\mathbb{A} \times_{\mathbb{F}_0 } \mathbb{F}_1 \rightarrow \mathbb{A},
\end{align*}    
the first one a left action and the second one a right action such that are {\em Courant actions}, meaning that their graphs are Dirac structures with support the graphs of the corresponding actions of the base groupoids. We say that these actions constitute a Courant Morita equivalence if they define a Morita equivalence in the usual sense \cite{dufzun} and $-p^*_0 \omega_0+ q^*_0 \omega_0'=dH$. 

In this work we shall only care about the following particular case of this concept. We shall say that a Courant Morita equivalence is {\em exact} if the twisted CA-groupoids are exact, the bimodule is of the form $\mathbb{A}=\mathbb{T}_H S$ and the Courant actions are given by exact Dirac structures \cite[\S 3]{baigua}. \end{defi} 

\begin{rema}\label{rem:mor} There is a useful way of packaging the information of a Courant Morita equivalence. The fibred product
\[ \mathbb{E}_1 \times_{\mathbb{E}_0 } \mathbb{A} \times_{\mathbb{F}_0 } \mathbb{F}_1 :=\{(e,a,f)\in \mathbb{E}_1 \times \mathbb{A} \times \mathbb{F}_1:\mathtt{s}(e)=p(a), \mathtt{t}(f)=q(a) \}      \] 
is a VB-groupoid $\mathbb{E}_1 \times_{\mathbb{E}_0 } \mathbb{A} \times_{\mathbb{F}_0 } \mathbb{F}_1 \rightrightarrows \mathbb{A} $ with $(g,s,h)\mapsto s$ as the source map, $(g,s,h)\mapsto g\cdot s\cdot h$ as the target and the multiplication is
\[ \mathtt{m}( (g,g'\cdot s \cdot h',h),(g',s,h'))=(\mathtt{m}(g,g'),s,\mathtt{m}(h',h)).   \] 
This VB-groupoid is embedded by means of $(e,a,f)\mapsto (e,f^{-1},e\cdot a\cdot f,a)$ as a twisted multiplicative Dirac structure inside the twisted CA-groupoid given by the product
\[ (\mathbb{E}_1 \times \overline{\mathbb{F}_1}  \times  \overline{ \mathbb{A} } \times \mathbb{A}\rightrightarrows \mathbb{E}_0 \times \mathbb{F}_0 \times \mathbb{A})=(\mathbb{E}_1 \times \overline{\mathbb{F}_1}\rightrightarrows \mathbb{E}_0 \times \mathbb{F}_0 )\times (\overline{ \mathbb{A} } \times \mathbb{A} \rightrightarrows \times \mathbb{A})\] 
with support the image of the morphism
\begin{align*}  &(G_1 \times_{G_0} S \times_{H_0} H_1 \rightrightarrows S) \hookrightarrow (G_1 \times H_1 \times S \times S \rightrightarrows G_0 \times H_0 \times S),\\ 
&(g,s,h)\mapsto(g,h^{-1},g\cdot s\cdot h,s). \end{align*}       
\end{rema} 
\begin{thm}\label{thm:cla} The correspondence established in Proposition \ref{thm:excagpd} induces a bijection between $H^2(\mathcal{A}^2(G_\bullet))$ and the exact Courant Morita equivalence classes of exact twisted CA-groupoids over $G_1 \rightrightarrows G_0$, which restrict to the trivial Morita equivalence on the base groupoids.  \end{thm}  
\begin{rema} Let us point out that the inclusion of the normalised subcomplex of the Bott-Shulman-Stasheff complex into the total complex is a quasi-isomorphism (see \cite[Rem. 3]{dupder}) so the space $H^2(\mathcal{A}^2(G_\bullet))$ has an intrinsic topological meaning depending only on the Morita equivalence class of $G_1 \rightrightarrows G_0$ \cite{behcoh}. Also, the condition on the base bimodule being isomorphic to $G_1$ itself may be removed but at the cost of introducing a Picard group action on $H^2(\mathcal{A}^2(G_\bullet))$. \end{rema}
\begin{proof}[Proof of Theorem \ref{thm:cla}] Suppose that $\Omega= \sum_i \omega_i $ and $\Omega'=\sum_i \omega_i'$ are 2-shifted 2-forms on a Lie groupoid $G_1 \rightrightarrows G_0$ which are cohomologous. By Proposition \ref{thm:excagpd}, we can assume that $\Omega'- \Omega=D \eta_0$ for $\eta_0\in \Omega^3(G_0)$. Let $\mathbb{E}_1 \rightrightarrows \mathbb{E}_0$ and $\mathbb{F}_1 \rightrightarrows \mathbb{F}_0$ be the CA-groupoids associated to ${\Omega } $ and $\Omega'$ respectively. We can define the Courant bimodule as $\mathbb{A}:=  \mathbb{T}_{\omega_1'+\mathtt{t}^*\eta_0 }G_1$ and the moment maps as $\mathtt{t}_{\mathbb{E}_1 }:\mathbb{A} \rightarrow \mathbb{E}_0 $, $\mathtt{s}_{\mathbb{F}_1 }:\mathbb{A} \rightarrow \mathbb{F}_0 $. Then the multiplication maps:
\begin{align*} &\mathtt{m}_{\mathbb{E}_1 }: \mathbb{E}_1 \times_{\mathbb{E}_0 } \mathbb{A} \rightarrow \mathbb{A} \\
& \mathtt{m}_{\mathbb{F}_1 }: \mathbb{A} \times_{\mathbb{F}_0 } \mathbb{F}_1 \rightarrow \mathbb{A}
    \end{align*} 
are Courant actions are desired. This is just a variation of an argument described in \cite{gaudir}.

Conversely, suppose that we have an exact Courant Morita equivalence between exact twisted CA-groupoids $\mathbb{E}_1 \rightrightarrows \mathbb{E}_0$ and $\mathbb{F}_1 \rightrightarrows \mathbb{F}_0$ over $G_1 \rightrightarrows G_0$ with bimodule $\mathbb{A} $ over $S$. Using splittings of $\mathbb{E}_1$ and $\mathbb{F}_1$ as in Proposition \ref{thm:excagpd}, we have that they correspond to 2-shifted 2-forms $\Omega= \sum_i \omega_i$ and $\Omega'= \sum_i \omega_i'$. Since the Courant actions are exact, they are given by the graphs of 2-forms and so their associated multiplicative Dirac structure defined by Remark \ref{rem:mor} is exact as well and so there is a normalised 2-form $\eta_1\in \Omega^2( G_1 \times_{G_0} S \times_{G_0} G_1)$ such that
\[ -d\eta_1=\text{pr}_1^* \omega_1 - \text{pr}_3^* \omega_1'- \mathtt{t}^* H +\mathtt{s}^*H.    \] 
If we put $\eta_0=H$, we have that $D (H-\eta_1)=-\text{pr}_1^* \Omega + \text{pr}_3^* \Omega'$ on the Lie groupoid $\mathcal{G} :={G_1} \times_{G_0} S \times_{G_0} G_1 \rightrightarrows S$. Since both projection maps from $\mathcal{G} $ to $G_1$ are Morita equivalences \cite[Prop. 7.2.10]{dufzun} and $S$ is assumed to be isomorphic to $G_1$ as bimodule, this Morita equivalence determines the identity isomorphism on the cohomology of the complex of 2-forms on $G_1\rightrightarrows G_0$, see \cite{behcoh,getdia}. It follows that $\Omega $ and $\Omega'$ have to be cohomologous as well. \end{proof}   
The cohomology space $H_{dR}^2(G_0)$ on a manifold $G_0$ classifies central Lie algebroid extensions of $TG_0$ by $\mathbb{R} $, the cohomology space $H^1(\mathcal{A}^2(G_\bullet))$ of 1-shifted 2-forms on a Lie groupoid $G_1 \rightrightarrows G_0$ classifies what we can call {\em quasi} (or twisted) LA-groupoid central extensions of $TG_1 \rightrightarrows TG_0$ by $\mathbb{R}$, regarded as a groupoid $\mathbb{R} \rightrightarrows \{0\}$.  The integral classes inside these cohomology groups appear to be related to prequantisations of 0 and 1-shifted symplectic groupoids \cite{bunger} so our result could correspond to a description of the Atiyah algebroids corresponding to prequantisations of 2-shifted (pre)symplectic groupoids \cite{getdia}.    

\section{The first Pontryagin class and twisted CA-groupoids over gauge groupoids}
\subsection{The foliation by Courant algebroids of a CA-groupoid}
Let $\mathbb{E}_1 \rightrightarrows \mathbb{E}_0$ be a not necessarily exact CA-groupoid over $G_1 \rightrightarrows G_0$ twisted by $\omega_0\in \Omega^4_{cl}(G_0)$. We can see this kind of object as a graded geometric analogue of a symplectic groupoid, see \cite{lacou}. Then we can expect that there is an induced symplectic foliation on its base in which the symplectic leaves are given by (twisted) Courant algebroids \cite{roycou}. This is made precise by the following result which is the global analogue of an observation in \cite{coutyp}.
\begin{thm}\label{thm:folca} Let $i:\mathcal{O} \hookrightarrow G_0$ be a $G_1$-orbit and take $p\in \mathcal{O} $. Then the source-fibre $\mathtt{s}_{\mathbb{E}_1 }^{-1}(0_p) \hookrightarrow \mathbb{E}_1 \rightrightarrows \mathbb{E}_0$ reduces to a $-i^* \omega_0$-twisted Courant algebroid over $\mathcal{O} $ given by the quotient of the vector bundle
\[ E_p:=\mathtt{s}_{\mathbb{E}_1 }^{-1}(0_p)/(\ker \mathtt{s}_{\mathbb{E}_1 }\cap \ker \mathtt{t}_{\mathbb{E}_1 })|_{\mathtt{s}^{-1}(p)}  \]
by the action of the isotropy group $G_1(p) =\mathtt{s}_{G_1}^{-1}(p)\cap \mathtt{t}_{G_1}^{-1}(p) \hookrightarrow G_1$ defined by right multiplication:
\[  v \cdot g := \mathtt{m}_{\mathbb{E}_1 }(v,0_g), \quad \forall v\in E_p, \, \forall g\in G_1(p).  \] 
\end{thm}      
\begin{proof} The fact that the graph of the multiplication is lagrangian implies that the kernel of the metric restricted to $\mathtt{s}_{\mathbb{E}_1 }^{-1}(0_p)$ is given exactly by $\ker \mathtt{s}_{\mathbb{E}_1 }\cap \ker \mathtt{t}_{\mathbb{E}_1 }|_{\mathtt{s}^{-1}(p)}$ and also that the metric descends to $E_p/G_1(p)$ as a well defined nondegenerate symmetric bilinear form. The involutivity of the graph of the multiplication immediately implies that the Courant bracket on $\mathbb{E}_1$ induces a Courant bracket on $E_p/G_1(p)$ which is $-i^* \omega_0$-twisted. \end{proof}
\begin{rema} The Courant algebroids obtained by reduction according to Theorem \ref{thm:folca} can be seen as the symplectic leaves of a degree two Poisson NQ-manifold as in \cite{coutyp,degtwo}. It would be interesting to see in which sense such a graded Poisson structure is obtained by reduction from $\mathbb{E}_1$. \end{rema}  
\begin{exa} In the case of the canonical Courant algebroid over a Lie groupoid, the Courant algebroids induced on the orbits as in the Theorem above are the corresponding canonical Courant algebroids. \end{exa} 
\subsection{Transitive Courant algebroids and exact CA-groupoids}\label{sec:tracou} Now we shall apply the observations that we made to the construction of heterotic Courant algebroids. 

Let $G$ be a Lie group with quadratic Lie algebra $\mathfrak{g} $ with $\langle \,,\, \rangle $ being its Ad-invariant symmetric non-degenerate pairing. If $P$ is a principal $G$-bundle, the kernel of the anchor of the Atiyah algebroid $A$ of $P$ is isomorphic to the adjoint bundle $\mathfrak{g} (P)=(P \times \mathfrak{g})/G $ and it inherits a metric from $\langle \,,\, \rangle $.
\begin{defi}[\cite{barhek}] Let $E$ be a Courant algebroid over a manifold $M$ such that its anchor $\mathtt{a} $ is surjective and such that the quotient Lie algebroid $A:=E/\mathtt{a}^*(T^*M)$ is isomorphic to the Atiyah algebroid of a principal $G$-bundle $\pi:P \rightarrow M$ by an isomorphism that preserves the restriction of the metrics to the kernel. Then we say that $E $ is a {\em heterotic Courant algebroid} and $E$ is called a Courant algebroid extension of $A$. \end{defi} 
Let us consider the closed 2-shifted 2-form \eqref{eq:bg} for $K=G$. A principal $G$-bundle $P$ over a manifold $M$ induces a Lie groupoid morphism from the action groupoid associated to $P$:
\[ \text{pr}_2: (P \rtimes G \rightrightarrows P) \rightarrow G; \]
since $H^\bullet_{dR}(M)$ is isomorphic to the Bott-Shulman-Stasheff cohomology of the groupoid $P\rtimes G \rightrightarrows P$, we can define the first Pontryagin class of $P$, denoted by $p_1(P)\in H^4_{dR}(M)$ as the class of the pullback $\text{pr}_2^* \Omega $. If one takes a principal connection $\nabla$ on $P$, then $p_1(P)=[\langle F_\nabla \wedge F_\nabla \rangle ]$, where $F_\nabla$ is the curvature of $\nabla$, see \cite[\S 4.2]{walmulger}. 

Now we will see how the vanishing of this cohomology class implies that the Atiyah algebroid $A$ of $P$ admits a heterotic Courant algebroid extension.
\begin{enumerate} \item Let $(P \times P)/G \rightrightarrows M$ be the gauge groupoid of $P$. Since $(P \times P)/G \rightrightarrows M$ is Morita equivalent to $G$, the class of $\Omega $ can be represented by a closed 2-shifted 2-form $\Omega'=\sum \omega_i'$. Proposition \ref{thm:excagpd} implies that $\Omega'$ corresponds to an exact $\omega'_0$-twisted CA-groupoid $\mathbb{E}_1 \rightrightarrows \mathbb{E}_0$ over $(P \times P)/G \rightrightarrows M$. 
\item Theorem \ref{thm:folca} implies then that there is a $-\omega'_0$-twisted Courant algebroid on $M$ which is the only $(P \times P)/G$-orbit and is a twisted Courant algebroid extension of $A$ by $T^*M$. We can see that $\omega_0'$ represents the first Pontryagin class of $P$ and hence we obtain part of our desired result.
\item Conversely, we can see that given a (twisted) Courant algebroid extension of $A$, there is a simple canonical construction of a (twisted) CA-groupoid over $\mathcal{G} $.
\end{enumerate}

\begin{thm}\label{thm:hetca} Let $\pi:P \rightarrow M$ be a principal $G$-bundle such that $\mathfrak{g} $ is a quadratic Lie algebra with pairing $\langle \,,\, \rangle $. Then the following statements hold.
\begin{enumerate} \item There is an exact CA-groupoid over the gauge groupoid $(P \times P)/G \rightrightarrows M$ that induces a transitive twisted Courant algebroid over $M$ as in Theorem \ref{thm:folca} with twisting class given by $-p_1(P)$. 
\item The Atiyah algebroid of $P$ admits a heterotic Courant algebroid extension $E $ if and only if $p_1(P)$ vanishes.
\item Every (twisted) heterotic Courant algebroid extension of $A$ is obtained by applying Theorem \ref{thm:folca} to an exact (twisted) CA-groupoid as in item (1). \end{enumerate} 
 \end{thm}
\begin{rema} Item (2) was originally discovered in \cite{sevlet}. Since then, several explanations have been given for this result, see \cite{brepon,pymsaf} and their references. Items (1) and (3) seem to be new and put this well known result into a more Lie theoretical context. \end{rema}  
\begin{proof}[Proof of Theorem \ref{thm:hetca}] {\em Step 1: the 2-shifted closed 2-form $\Omega $ on $G_\bullet$ defined in \eqref{eq:bg} induces an exact twisted CA-groupoid over $(P \times P) /G \rightrightarrows M$}. This is a consequence of the fact that any inclusion of $i_p:G \hookrightarrow (P \times P)/G $ as an isotropy group over $p\in M$ is a Morita equivalence and Morita equivalences induce isomorphisms in the Bott-Shulman-Stasheff cohomology \cite{behcoh,getdia}. So the 2-shifted 2-form \eqref{eq:bg} can be transported to a closed 2-shifted 2-form $\Omega'=\sum \omega_i'$ such that $i_p^*[\Omega']=[\Omega ]$. Since the normalised subcomplex of the Bott-Shulman-Stasheff complex is homotopy equivalent to the total complex \cite[Rem. 3]{dupder}, we can assume that $\Omega'$ is normalised. Notice that, based on Theorem \ref{thm:cla}, this fact is enough to show that item (1) holds. 

For the sake of clarifying why (2) holds as well, let us sketch how $i_p$ induces an isomorphism $\mathcal{A}^2(i_p)$ between the cohomology of the complexes of 2-forms $\mathcal{A}^2 (G_\bullet)$ and $\mathcal{A}^2( ((P \times P)/G)_\bullet)$ \cite{behcoh}. Consider the Lie groupoid morphism $\Phi: P^2 \times G \rightarrow P^2/G$ from the product of the pair groupoid with $G$ defined by $(p,q,a)\rightarrow [pa,q]$ for all $(p,q,a)\in P \times P \times G$; $\Phi$ is a Morita morphism and so there is $\Omega'\in \mathcal{A}^2( ((P \times P)/G)_\bullet)$ such that $[\Phi^* \Omega']= [\text{pr}_2^*\Omega  ]$, where $\text{pr}_2: P^2 \times G \rightarrow G$ is the projection. Considering the pullback of $\text{pr}_2^*\Omega$ to the subgroupoid 
\[ \{ (x,xa,a)\in P \times P \times G|x\in P,\, a\in G\} \hookrightarrow P^2 \times G \]
which is isomorphic to the action groupoid $P \rtimes G \rightrightarrows P$, we see that $[\Phi^* \omega_0' ]\in H^4(\text{Tot} (\Omega^\bullet((P \times G))_\bullet,D)\cong H^4_{dR}(M)$ represents $p_1(P)$ by definition since it corresponds to the pullback of $\Omega$ along the projection $P \rtimes G \rightarrow G$. As a consequence of Theorem \ref{thm:cla}, we get that there is a representative of $[\Omega']$ that admits an untwisted CA-groupoid realization if and only if $p_1(P)=0$.

{\em Step 2: the twisted CA-groupoid $\mathbb{E}_1 \rightrightarrows \mathbb{E}_0$ induced by $\Omega'$ satisfies that $\ker \mathtt{s}_{\mathbb{E}_1 }\cap \ker\mathtt{t}_{\mathbb{E}_1 }=0$}. This is because the source and target maps on $\mathbb{E}_1=\mathbb{T}_{\omega_1'}((P \times P)/G) $ are defined by \eqref{eq:prstmap} and then, over $p\in M$, an element $X\oplus \alpha \in(\ker \mathtt{s}_{\mathbb{E}_1 })_p\cap (\ker\mathtt{t}_{\mathbb{E}_1 })_p$ satisfies that 
\begin{align} \langle s_0^*i_{(X,0 ) }\omega_2',u \rangle+\langle s_1^*i_{(0,X ) }\omega_2',v \rangle=-\langle \mathtt{s}_{T^*G} (\alpha)  ,u \rangle -\langle \mathtt{t}_{T^*G} (\alpha)  ,v \rangle\quad \forall u\in \ker T_p\mathtt{t},\, v\in \ker T_p\mathtt{s}. \label{eq:nondegpai}   
\end{align}  
Then for $u=-v=Y\in \ker T_p \mathtt{s}\cap \ker T_p \mathtt{t}$ we have that
\[ \omega_2'((X,0),(0,Y))-\omega_2'((0,X),(Y,0))=\langle \alpha,Y \rangle - \langle \alpha , Y \rangle=0. \]
On the other hand, $i_p^*[\Omega']=[\Omega ]$ and the pairing defined above is independent of the choice of representative in the cohomology class since it is given by the component $\lambda^{\Omega'}_1$ of the chain complex morphism induced by $\Omega'$ as in \eqref{eq:getpai}, see \cite[Lemma E.1]{cuezhu}. Since the pairing induced as above by $\Omega $ is the bilinear form $-\langle \,,\, \rangle $ \cite[Prop. 3.5]{cuezhu}, we have that it is nondegenerate when restricted to the isotropy Lie algebras. So in our case $X=0$ and hence \eqref{eq:nondegpai} applied to $v=u^{-1}=-u+\mathtt{a}(u)$ implies that $\alpha =0$ as well.

{\em Step 3: the twisted Courant algebroid over $M$ induced by Theorem \ref{thm:folca} is transitive and it is an extension of the Lie algebroid $A$ of $(P \times P)/G \rightrightarrows M$.} Since $\ker \mathtt{s}_{\mathbb{E}_1 }\cap \ker\mathtt{t}_{\mathbb{E}_1 }=0$, Theorem \ref{thm:folca} implies that the quotient 

\[ E:= \mathtt{s}_{\mathbb{E}_1 }^{-1}(0_x)/G \]
for any $x\in M$ is an $-\omega_0'$-twisted Courant algebroid which is an extension of $A$. In fact, the anchor map of $\mathbb{E}_1$ induces a surjective bracket preserving map $E \rightarrow A$ as desired. Therefore, our claim holds.  

{\em Step 4: there is a canonical exact action CA-groupoid over the gauge groupoid of a heterotic Courant algebroid}. Let $E $ be a heterotic Courant algebroid associated to $\pi: P \rightarrow M$. Then the VB-groupoid defined by
\begin{align} \mathtt{t}^*E\oplus \mathtt{s}^* E\rightrightarrows E   \label{eq:exca} \end{align} 
is an exact CA-groupoid over $(P \times P)/G \rightrightarrows M$ if we equip it with the metric of $E$ on the first summand and its negative on the second one. The Courant bracket and the anchor are obtained by observing that $E \times \overline{E }  $ acts on the map $(\mathtt{t},\mathtt{s}):(P \times P)/G \rightarrow M \times M  $ by means of $(u,v) \mapsto \rho(u)^r+(\rho(v)^{-1})^l$ for all $(u,v) \in \Gamma (E )\oplus \Gamma (  E)$, where $\rho:E \rightarrow A$ is the quotient map. Since this action has lagrangian stabilizers, given at each point $g$ by $\{\xi\oplus \eta\in E_{\mathtt{t}(g) }\oplus E_{\mathtt{s}(g) }  | \rho(\xi)^r+(\rho(\eta)^{-1})^l=0\} $, we have that \cite[Thm. 2.12]{coupoi} implies the result and by dimension counting we get that this Courant algebroid is exact. It is immediate to see that the leafwise Courant algebroid induced by this twisted CA-groupoid as in Theorem \ref{thm:folca} is $E$ again. By splitting $\mathtt{t}^*E\oplus \mathtt{s}^* E \rightrightarrows E$ we can obtain an explicit 2-shifted 2-form on $(P \times P)/G \rightrightarrows M$ that classifies it as in Theorem \ref{thm:cla}. \end{proof} 

\begin{rema} We can interpret the previous result as an integration of a (twisted) heterotic Courant algebroid to a 2-shifted symplectic structure $\Omega'$ \cite{getdia,bunger} on the associated gauge groupoid (in accordance with the observations in \cite{pymsaf}), see Definition \ref{def:2shisym}. In order to make this interpretation precise, we have to apply the Van Est map to $\Omega' $ as in \cite{cuezhu}. The nondegeneracy of $\Omega'$ is encoded in the fact that the associated twisted CA-groupoid has trivial isotropy as in Step 2 of the previous proof. If $M$ is a point, then $E=A=\mathfrak{g}  $ and what we recover is the description of an action Courant algebroid on $G$ that corresponds to \eqref{eq:bg} by means of a splitting \cite{purspi}. \end{rema} 

\begin{rema} Let us compare our approach with the construction of \cite[Prop. 3.5]{barhek} and \cite[Thm. 5]{tduacou}. Let us continue using the notation of the proof of Theorem \ref{thm:hetca}. We will see that indeed any heterotic Courant algebroid can be obtained by reduction of an exact Courant algebroid over the associated principal bundle. According to the discussion in Step 4, we can represent the CA-groupoid $\mathbb{E}_1 \rightrightarrows \mathbb{E}_0$ associated to $\Omega'$ as in \eqref{eq:exca}. Then the pullback Courant algebroid $\mathbb{E}_1|_P$ of $\mathbb{E}_1$ \cite{coupoi} to the source fibre $P\cong \mathtt{s}^{-1}(x) \hookrightarrow (P \times P)/G$ is exact and it is given by
\[ \mathbb{E}_1|_P=\{\xi\oplus \eta\in E_{\mathtt{t}(g) }\oplus E_{\mathtt{s}(x) }|\rho(\eta)\in \ker \mathtt{a}_A ,\, g=[p,x]\in \mathtt{s}^{-1}(x)     \}; \]
where $\rho:E \rightarrow A$ is the quotient map. We have that $\mathbb{E}_1|_P\cong \mathtt{t}^* E\oplus \mathfrak{g}  $ and, as we know from the the proof of Theorem \ref{thm:folca}, $\ker \mathtt{s}_{\mathbb{E}_1 }$ is the vector bundle that reduces to $E $ by the $G$-action. So in order to reduce $\mathbb{E}_1|_P$ we just to need to observe that the canonical map $\psi:\mathfrak{g} \rightarrow \mathbb{E}_1|_P $ is an extended $\mathfrak{g} $-action and by taking $\psi(\mathfrak{g} )^\perp/G$ we get $E $ again. \end{rema} 

\begin{exa} Let $D$ be a Lie group with quadratic Lie algebra $\mathfrak{d} $ and take $C \hookrightarrow D$ a closed subgroup such that $\mathfrak{c} \hookrightarrow \mathfrak{d} $ is coisotropic and $\mathfrak{c}^\perp \hookrightarrow \mathfrak{c}$ integrates to a closed connected subgroup $C^\perp \hookrightarrow C$. The canonical action of $\mathfrak{d} $ on the homogeneous space $M=D/C$ (defined by the right action of $C$ on $D$) has coisotropic stabilizers and so there is a transitive Courant algebroid structure on $E:= \mathfrak{d} \times M$ \cite{coupoi}. Consider the principal $C/C^\perp$-bundle $P=D/C^\perp \rightarrow M$. So we see that $C/C^\perp$ can play the role of $G$ in Theorem \ref{thm:hetca}. The Atiyah algebroid of $P$ is $A= T (D/C^\perp)/(C/C^\perp)$ so the map $\rho: E \rightarrow A$, $\rho(u,xC)=[\left.\frac{d}{dt}\right|_{t=0}e^{tu}xC^\perp ]$ is the desired quotient map, this is well defined because $C^\perp \hookrightarrow C$ is a normal subgroup. The gauge groupoid associated to $P$ is isomorphic to a homogeneous space as well:
\[ (D \times D) /L \rightrightarrows D/C, \quad L  =\{(x,ax)\in D \times D |x\in C,a\in C^\perp \} ;\]   
where the $L$-action on $D \times D$ is by right translation ($L \hookrightarrow C$ is a Lie 2-subgroup of the pair 2-group $D \times D \rightrightarrows D$). The exact CA-groupoid $\mathtt{t}^*E\oplus \mathtt{s}^* E$ that we get on top of the homogeneous space $(D \times D) /L$ is then also an action Courant algebroid, given by the canonical action of $\mathfrak{d} \times \overline{\mathfrak{d}}  $ on $(D \times D)/L$ which has coisotropic (lagrangian) stabilizers. \end{exa} 
\section{Dirac structures and 2-shifted lagrangian groupoids} In this section we shall introduce the double quasi-symplectic groupoids which serve as integrations of multiplicative Dirac structures in exact CA-groupoids. We shall illustrate this concept with some examples coming from Dirac structures in transitive Courant algebroids. Finally, we shall see that double quasi-symplectic groupoids induce local 2-shifted symplectic groupoids.
\subsection{Multiplicative Dirac structures and double quasi-symplectic groupoids} Let $L_1 \rightrightarrows L_0$ be a multiplicative Dirac structure with full support inside the exact (untwisted) CA-groupoid over $G_1 \rightrightarrows G_0$ with associated closed 2-shifted 2-form given by $\omega =\omega_2+\omega_1$. It is natural to ask what is the global counterpart of such a multiplicative Dirac structure.

Consider a groupoid object in the category of Lie groupoids, such a structure is denoted by a diagram of the following kind
\[ \xymatrix{ X_{1,1} \ar@<-.5ex>[r] \ar@<.5ex>[r]\ar@<-.5ex>[d] \ar@<.5ex>[d]& X_{1,0}  \ar@<-.5ex>[d] \ar@<.5ex>[d] \\ X_{0,1}  \ar@<-.5ex>[r] \ar@<.5ex>[r] & X_{0,0}, }\]  
where each of the sides represents a groupoid structure and the structure maps of $X_{1,1} $ over $X_{1,0}$ are groupoid morphisms with respect to $X_{1,1} \rightrightarrows X_{0,1}$ and $X_{1,0} \rightrightarrows X_{0,0}$. A {\em double Lie groupoid} is an object as above such that the map $(\mathtt{t}^V,\mathtt{s}^H):  X_{1,1} \rightarrow X_{0,1} \times_{X_{0,0}} X_{1,0} $ is a submersion (the superindices $\quad^H,\quad^V$ denote the groupoid structures $X_{1,1} \rightrightarrows X_{1,0}$, $X_{1,1} \rightrightarrows X_{0,1}$ respectively); if $(\mathtt{t}^V,\mathtt{s}^H)$ is in addition surjective, we call this structure a full double Lie groupoid \cite{browmac,macdou}. Out of a double Lie groupoid  we can construct its {\em nerve} which is the bisimplicial manifold $X_{a,b}$ obtained by taking the nerves of both vertical and horizontal groupoid structures: 
\[ X_{a,b}=\{(\xi_{i,j})|\xi_{i,j}\in X_{1,1},\, \mathtt{s}^V(\xi_{i-1,j})=\mathtt{t}^V(\xi_{i,j}),\, \mathtt{s}^H(\xi_{i,j-1})=\mathtt{t}^H(\xi_{i,j}),\, 1\leq i\leq a, \, 1\leq j\leq b\}; \]
for the sake of simplifying notation, we shall identify a double Lie groupoid with its nerve. 
\begin{defi}\label{def:douqua} Let $\mathfrak{X}=(X_{a,b})$ be a double Lie groupoid. By taking the two associated simplicial differentials on its nerve with respect to the horizontal and vertical groupoid structures $\delta^H,\delta^V$ we get a double complex $(C_{a,b}=\Omega^\bullet(X_{a,b} ),\delta^H,\delta^V)$. We say that $\mathfrak{X} $ is a {\em double quasi-symplectic groupoid} if there are normalised differential forms $\Omega_{1,1} \in \Omega^2(X_{1,1})$, $\Omega_{2,0} \in \Omega^2(X_{2,0})$ and $\Omega_{1,0}\in \Omega^3(X_{1,0})$ such that the following two conditions hold. 
 Consider the de Rham differential $d$ on $C_{a,b}$ and define the total differential 
\[ \widetilde{D}=\delta^V+(-1)^a\delta^H+(-1)^{a+b}d. \] 
So we have that $\Omega:=\Omega_{1,0}+\Omega_{1,1}+\Omega_{2,0}$ is a quasi-symplectic structure if (1) $\widetilde{D} \Omega =0$ and (2) it is horizontally nondegenerate:
\[ \ker \Omega_{1,1}|_x\cap \ker T\mathtt{t}^H|_x \cap\ker T \mathtt{s}^H|_x=0 \quad \forall x\in X_{1,0} ; \quad \dim X_{1,1}=2\dim X_{1,0}.  \] 
We denote such a double quasi-symplectic groupoid as $(\mathfrak{X},\Omega)$. \end{defi}
\begin{exa} A double quasi-symplectic groupoid $(\mathfrak{X},\Omega)$ is a double symplectic groupoid \cite{luwei2} if $\Omega=\Omega_{1,1} $ is symplectic. \end{exa}  
\begin{exa}[AMM groupoid] Let $K$ be a Lie group with a quadratic Lie algebra. The action groupoid $K\times K\rightrightarrows K$
\[ \mathtt{t}^H(a,x)=axa^{-1},\quad \mathtt{s}^H(a,x)=x, \quad \mathtt{m}^H((a,bxb^{-1}),(b,x))=(ab,x),   \] 
together with the Lie group bundle structure over $K$
\[ \mathtt{t}^V(a,x)=\mathtt{s}^V(a,x)=a, \quad \mathtt{m}^V((a,x),(a,y))=(a,xy);   \]
constitutes a double Lie groupoid:
\[ \xymatrix{ K \times K \ar@<-.5ex>[r] \ar@<.5ex>[r]\ar@<-.5ex>[d] \ar@<.5ex>[d]& K  \ar@<-.5ex>[d] \ar@<.5ex>[d] \\ K  \ar@<-.5ex>[r] \ar@<.5ex>[r] & \{\ast\}. }\]
This double Lie groupoid becomes quasi-symplectic with the forms
\[ \Omega_{1,1}|_{(k,x)}= \langle \mathrm{Ad}_x \text{pr}_1^*\theta^l\wedge\text{pr}_1^*\theta^l\rangle+\langle\text{pr}_1^*\theta^l\wedge\text{pr}_2^*(\theta^l+\theta^r)\rangle,\quad \Omega_{1,0}=\Omega_1, \quad \Omega_{2,0}=\Omega_2; \]
using the formulas of \eqref{eq:bg} and denoting by $\theta^l$ and $\theta^r$ the left and right invariant Maurer-Cartan forms on $K$. This double quasi-symplectic groupoid integrates the Cartan-Dirac structure \cite[\S 7.2]{burint}. \end{exa}
By applying the Lie functor to a double Lie groupoid in one direction we get an LA-groupoid: a groupoid in the category of Lie algebroids \cite{macdou}. In the case of a double quasi-symplectic groupoid, we get a multiplicative Dirac structure. We omit the unit embeddings in what follows, for simplicity.
\begin{prop}\label{pro:difdir} Take a double quasi-symplectic groupoid as in Definition \ref{def:douqua}, then its horizontal differentiation is a multiplicative Dirac structure inside the exact CA-groupoid over $X_{1,0} \rightrightarrows X_{0,0}$ defined by the closed 2-shifted 2-form $\Omega_{2,0}+\Omega_{1,0}$. \end{prop}
\begin{proof} Based on \cite{burint}, it only remains to show that the horizontal LA-groupoid $A^H \rightrightarrows A_{0,1}$ is embeddded by the map $u\mapsto F(u):=T\mathtt{t}^H(u)\oplus i_u\Omega_{1,1} $ as a VB-subgroupoid of $\mathbb{T}_{\Omega_{1,0}}X_{1,0}$ with respect to the twisted multiplication determined by $\Omega_{2,0}$ as in Proposition \ref{thm:excagpd}. Let us observe first that $ i_{T\mathtt{s}^V(u) }\Omega_{1,1}\in T^\circ X_{0,0}\cong A_{1,0}^* $ for all $u\in TX_{1,1}$ since $\Omega_{1,1}$ is vertically normalised, so $F(A_{0,1})\subset TX_{0,0}\oplus A_{1,0}^*$. Now take $u \in A^H$, by definition we have that 
\[ \langle p_{T^*}\mathtt{s}_{\Omega_{2,0}}F(u), w \rangle = \langle i_u \Omega_{1,1},T \mathtt{m}(0,w)   \rangle+\Omega_{2,0}((T \mathtt{t}^H (u),T\mathtt{s}_{X_{1,0}}T\mathtt{t}^H(u)),(0,w) )    \]
where $p_{T^*}$ is the projection to $T^*X_{1,0}$ and $w\in A_{1,0}\cong \ker T \mathtt{t}_{X_{1,0}}$. But from evaluating the two terms in the equation $\delta^V \Omega_{1,1}=-\delta^H \Omega_{2,0}$ on $(u,T\mathtt{s}^V(u))$ and $(0,w)$ we get that 
\[ \langle p_{T^*}\mathtt{s}_{\Omega_{2,0}}F(u), w \rangle= i_{T\mathtt{s}^V(u) }\Omega_{1,1}(w) \]  
as desired, so $\mathtt{s}_{\Omega_{2,0}}F(u)=F(T\mathtt{s}^V (u))$. Similarly, we can see that $F$ interwines the target maps. Finally, if $u,u'\in A^H$ and $w,w'\in TX_{1,0}$ are suitably composable, then our definition means that
\begin{align*} &\langle p_{T^*} \mathtt{m}_{\Omega_{2,0}}(F(u),F(u')) ,T\mathtt{m}_{X_{1,0}}(w,w') \rangle=\\
&=\langle i_u \Omega_{1,1}, w \rangle + \langle i_{u'} \Omega_{1,1},w' \rangle -\Omega_{2,0}((T \mathtt{t}^H(u),T \mathtt{t}^H(u')  ),(w,w')) \end{align*}  but from $\delta^V \Omega_{1,1}=-\delta^H \Omega_{2,0}$ we get that $F$ is a groupoid morphism as we wanted. \end{proof}
Notice that integrating multiplicative Dirac structures would require understanding the integration of LA-groupoids which is currently an open problem. 

\subsection{2-shifted lagrangian groupoids and quasi-symplectic groupoids} Now we will see how Dirac structures inside transitive Courant algebroids allow us to construct double quasi-symplectic groupoids and Morita equivalences between them.
 
Let $L \hookrightarrow E$ be a Dirac structure inside a heterotic Courant algebroid $E$ over $M$ with support an immersed submanifold $N \hookrightarrow M$. According to \cite{pymsaf}, $L$ should be seen as an infinitesimal 2-shifted lagrangian Lie algebroid inside $E$, regarded as the infinitesimal data of a 2-shifted symplectic groupoid. In this work we shall begin by formulating all the relevant concepts in a Courant-theoretic language, their interpretation in the framework of shifted symplectic geometry arises after introducing splittings. 
\begin{defi}\label{def:2lag} Let $\mathtt{t}^* E\oplus \mathtt{s}^* E \rightrightarrows E $ be the exact CA-groupoid over a gauge groupoid $(X_1 \rightrightarrows X_0)=((P \times P)/G \rightrightarrows M)$ associated to a heterotic Courant algebroid $E $ over $M$ and let $\mathcal{L}_1 \rightrightarrows X_0  $ be a Lie groupoid. We say that $\mathcal{L}_1 \rightrightarrows X_0  $ is a 2-shifted lagrangian groupoid inside $X_1 \rightrightarrows X_0$ if there is an exact morphism of CA-groupoids  
\[ R: \mathbb{T}\mathcal{L}_1 \rightarrow  \mathtt{t}^* E\oplus \mathtt{s}^*E  \] 
over some Lie groupoid morphism $\Phi: \mathcal{L}_1 \rightarrow X_1$ such that $\ker R\cap( T \mathcal{L}_1 \oplus 0 ) =0$ and $\dim \mathcal{L}_1=2\dim X_0 + \frac{1}{2}\dim G$.
\end{defi}
Let us observe that, by identifying $\mathtt{t}^* E\oplus \mathtt{s}^* E\rightrightarrows E $ as a CA-groupoid coming from a closed 2-shifted 2-form $(\mathbb{T}X_1 \rightrightarrows TX_0\oplus A_{X_1}, \omega =\omega_1+\omega_2) $, we can apply Corollary \ref{cor:exdir} to the previous definition. In fact, $R$ as in Definition \ref{def:2lag} is, in particular, an exact multiplicative Dirac structure $R\subset \overline{\mathtt{t}^* E\oplus \mathtt{s}^*E} \times \mathbb{T}\mathcal{L}_1  $ supported on the graph of $\Phi$. As a result, a 2-shifted lagrangian groupoid $\Phi:(\mathcal{L}_1 \rightrightarrows X_0) \hookrightarrow (X_1 \rightrightarrows X_0)$ is given by a normalised 2-form $\beta_1\in \Omega^2(\mathcal{L}_1)$ such that 
\begin{align}  &D \beta_1= \Phi^* \omega,  \notag\\
& \ker \beta_1\cap \ker T\Phi=0; \label{eq:2lag1}\\
&\dim \mathcal{L}_1=2\dim X_0 + \frac{1}{2}\dim G ; \notag  \end{align}
in this situation, $\beta_1$ is a {\em 2-shifted lagrangian structure on $\Phi$} \cite{ptvv}. Using this explicit description, one can identify a number of examples in Poisson geometry \cite{2lagpoi}; on the other hand, an advantage of Definition \ref{def:2lag} is that it does not require the choice of a splitting of the corresponding Courant algebroids. When $E$ is an exact Courant algebroid and $\omega=\omega_1$ is its classifying closed 3-form we recover the definition of quasi-symplectic (or pre-symplectic) groupoid \cite{mommapmor,burint}. We could generalise this definition to include twisted transitive Courant algebroids but we do not know any examples that would motivate this generality. 
\begin{prop}\label{pro:2lagdir} Consider a 2-shifted lagrangian groupoid $(\mathcal{L}_1 \rightrightarrows X_0,R,\Phi ) $ with respect to the exact CA-groupoid $\mathtt{t}^* E\oplus \mathtt{s}^* E \rightrightarrows E $ over $X_1 \rightrightarrows X_0$ associated to a heterotic Courant algebroid $E \rightarrow X_0$. Then the Lie algebroid $L$ of $\mathcal{L}_1 $ embeds canonically inside $E$ as a Dirac structure by means of the map $F:L \rightarrow E$ defined by 
\[ \text{graph}(F)= R\cap( (\mathtt{t}^* E\oplus 0)|_{X_0} \times   \ker T\mathtt{s}|_{\mathcal{L}_0 }  ).   \] \end{prop}  
\begin{proof} In terms of a splitting of $\mathtt{t}^* E\oplus \mathtt{s}^* E$ with $\omega=\omega_2+\omega_1$ as the associated closed 2-shifted 2-form and $\beta_1\in \Omega^2(\mathcal{L}_1)$ as the lagrangian structure, $R$ can be expressed as the Dirac structure
\[ R=\{(T\Phi(U)\oplus \xi,U\oplus \Phi^*\xi-i_U\beta_1)\in \mathbb{T}X_1\times \mathbb{T}\mathcal{L}_1 |U\oplus \xi\in T\mathcal{L}_1\oplus \Phi^*X_1\}\]
which is supported on $\text{Graph}(\Phi)$. Then $F$ as in the statement of the proposition is defined by 
\begin{equation}
  a\mapsto \phi(a)\oplus -s_0^*\left(i_{(\phi(a),0)}\omega_2\right)+  \mathtt{u}^*(i_a \beta)\in  \ker \mathtt{s}_\omega|_{X_0}\cong (\mathtt{t}^* E\oplus 0)|_{X_0}, \quad \forall a\in L;  \label{eq:embcou}
\end{equation}
where $\phi$ is the Lie algebroid morphism induced by $\Phi$ and we view the two terms in $-s_0^*\left(i_{(\phi(a),0)}\omega_2\right)+  \mathtt{u}^*(i_a \beta)$ as the respective components of a covector in $T^*X_1|_{X_0}= \text{Ann}(TX_0)\oplus \text{Ann}(\ker T\mathtt{t}|_{X_0})\cong \ker (T\mathtt{t}|_{X_0})^*\oplus T^*X_0$, so the result follows immediately from \eqref{eq:2lag1}, in particular, the involutivity of $F(L)$ follows from the fact that $R$ is a Courant relation. \end{proof} 
\begin{rema} As a consequence of Proposition \ref{pro:2lagdir}, we have then that $R$ is a Manin pair morphism from $(\mathbb{T} \mathcal{L}_1, T \mathcal{L}_1 )  $ to $(\mathtt{t}^* E\oplus \mathtt{s}^*E,\mathtt{t}^* L \oplus \mathtt{s}^*L)$ \cite{buriglsev}. \end{rema} 
Given two Dirac structures $L,L' \hookrightarrow E $, we can promote them to a single Dirac structure inside $\mathtt{t}^* E\oplus \mathtt{s}^* E$:
\[ \mathtt{t}^* L \oplus \mathtt{s}^* L' \hookrightarrow \mathtt{t}^* E\oplus \mathtt{s}^* E.  \]
It is natural to ask whether we can construct an integration for $\mathtt{t}^* L \oplus \mathtt{s}^* L' $ having some integrations of $L, L'$ at our disposal. It turns out that we can produce such an integration as an action Lie groupoid. Let $\mathcal{L}_1 \rightrightarrows X_0$ and $\mathcal{L}_1' \rightrightarrows X_0$ be 2-shifted lagrangian groupoids in the gauge groupoid $X_1 \rightrightarrows X_0$ associated to a heterotic Courant algebroid $E$ with moment maps $\Phi:\mathcal{L}_1 \rightarrow X_1$ and $\Phi':\mathcal{L}'_1 \rightarrow X_1$. Let us define the Lie groupoid $\mathcal{L}_1 \times_{X_0} X_1 \times_{X_0} \mathcal{L}_1' \rightrightarrows X_1$, where
\[ \mathcal{L}_1 \times_{X_0} X_1 \times_{X_0} \mathcal{L}_1'=\{(u,y,v)\in \mathcal{L}_1 \times X_1 \times \mathcal{L}_1'| \mathtt{s}(u)=\mathtt{t}(y),\, \mathtt{s}(y)=\mathtt{s}(v)    \}. \]
We think of this groupoid as the space of triples that induce a commutative square in 
$X_1$ so its groupoid structure is the action groupoid structure given by the projection $(u,y,v)\mapsto y$ as the source, $(u,y,v)\mapsto \Phi(u)y\Phi'(v^{-1})$ as the target and the multiplication is:
\begin{align}  \mathtt{m} ((u',\Phi(u)x\Phi'(v^{-1}),v') ,(u,x,v))=( u'u,x,v'v). \label{eq:douquamul}\end{align} 
In this situation, we can consider then three Dirac structures on $X_1$:
\[ \mathtt{t}^* L \oplus \mathtt{s}^* L, \quad \mathtt{t}^* L \oplus \mathtt{s}^* L', \quad \mathtt{t}^* L' \oplus \mathtt{s}^* L'; \] 
and using $\mathcal{L}_1,\mathcal{L}_1'$ we can produce integrations for all of them that are related as follows. We will assume that a splitting of $\mathtt{t}^* E \oplus \mathtt{s}^* E $ is chosen so that there is an associated closed 2-shifted 2-form $\omega =\omega_1+ \omega_2$ on $X_1 \rightrightarrows X_0$ and the 2-shifted lagrangian structures are encoded in $\beta_1\in \Omega^2(\mathcal{L}_1)$, $\beta_1'\in \Omega^2(\mathcal{L}_1')$ such that \eqref{eq:2lag1} holds.   
\begin{prop}\label{thm:2lagint} Let $(\mathcal{L}_1 \rightrightarrows X_0,\Phi,\beta_1)$ and $(\mathcal{L}_1' \rightrightarrows X_0,\Phi',\beta_1')$ be 2-shifted lagrangian groupoids in the gauge groupoid $X_1 \rightrightarrows X_0$ associated to a heterotic Courant algebroid $E $. Then the following statements hold.
\begin{enumerate} \item There are natural Lie groupoid structures $\mathcal{L}_1 \times_{X_0} X_1 \times_{X_0} \mathcal{L}_1 \rightrightarrows \mathcal{L}_1$, $\mathcal{L}_1' \times_{X_0} X_1 \times_{X_0} \mathcal{L}_1' \rightrightarrows \mathcal{L}_1'$ that together with the groupoid structures defined by \eqref{eq:douquamul} and the forms
\begin{align*}  &\Omega_{1,1}|_{(u,y,v)}=u^*\beta_1-v^*\beta_1+{(x,\Phi(v))}^*\omega_2-{(\Phi(u),y)}^*\omega_2, \\ 
&\Omega_{2,0}=\omega_2, \quad \Omega_{1,0}=\omega_1, \quad x\Phi(v)=\Phi(u)y, \\
&\Omega_{1,1}'|_{(u',y,v')}={u'}^*\beta_1'-{v'}^*\beta_1'+{(x,\Phi'(v'))}^*\omega_2-{(\Phi'(u'),y)}^*\omega_2, \\ 
&\Omega_{2,0}'=\omega_2, \quad \Omega_{1,0}'=\omega_1, \quad x\Phi'(v')=\Phi'(u')y, \\
&\forall (u,y,v)\in \mathcal{L}_1 \times_{X_0} X_1 \times_{X_0} \mathcal{L}_1,\quad \forall (u',y,v')\in \mathcal{L}_1' \times_{X_0} X_1 \times_{X_0} \mathcal{L}_1';
 \end{align*}  
make them into double quasi-symplectic groupoids that integrate the multiplicative Dirac structures $\mathtt{t}^* L \oplus \mathtt{s}^* L$ and $\mathtt{t}^* L' \oplus \mathtt{s}^* L'$ respectively.
\item The double quasi-symplectic groupoids of the previous item fit into a Morita equivalence of double Lie groupoids:

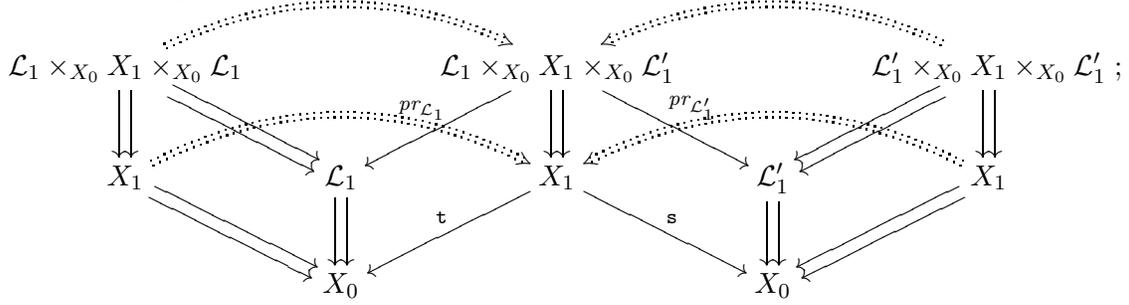
\[ \xymatrix{ \mathcal{L}_1 \times_{X_0} X_1 \times_{X_0} \mathcal{L}_1 \ar@<-.5ex>[rd] \ar@<.5ex>[rd] \ar@<-.5ex>[d] \ar@<.5ex>[d] \ar@{:>}@/^2pc/[rr] && \mathcal{L}_1 \times_{X_0} X_1 \times_{X_0} \mathcal{L}_1' \ar@<-.5ex>[d] \ar@<.5ex>[d] \ar[rd]^-{\text{pr}_{\mathcal{L}_1' } } \ar[ld]_-{\text{pr}_{\mathcal{L}_1 } }  && \mathcal{L}_1' \times_{X_0} X_1 \times_{X_0} \mathcal{L}_1' \ar@<-.5ex>[d] \ar@<.5ex>[d] \ar@<-.5ex>[ld] \ar@<.5ex>[ld] \ar@{:>}@/_2pc/[ll] \\
X_1 \ar@<-.5ex>[rd] \ar@<.5ex>[rd] \ar@{:>}@/^2pc/[rr] & \mathcal{L}_1 \ar@<-.5ex>[d] \ar@<.5ex>[d] & X_1\ar[rd]^-{\mathtt{s} } \ar[ld]_-{\mathtt{t} } & \mathcal{L}_1' \ar@<-.5ex>[d] \ar@<.5ex>[d] &X_1 \ar@<-.5ex>[ld] \ar@<.5ex>[ld] \ar@{:>}@/_2pc/[ll] \\
& X_0 &  & X_0 }; \]
where the dotted arrows represent the Lie groupoid actions defined by:
\begin{align*}  &\alpha^R: (\mathcal{L}_1 \times_{X_0} X_1 \times_{X_0} \mathcal{L}_1') \times_{\mathcal{L}_1' } (\mathcal{L}_1' \times_{X_0} X_1 \times_{X_0} \mathcal{L}_1') \rightarrow \mathcal{L}_1 \times_{X_0} X_1 \times_{X_0} \mathcal{L}_1', \\
&\alpha^R((u,y,v) ,(v,y',w))= (u, yy',w), \\ 
&\forall ((u,y,v) ,(v,y',w))\in (\mathcal{L}_1 \times_{X_0} X_1 \times_{X_0} \mathcal{L}_1') \times_{\mathcal{L}_1' } (\mathcal{L}_1' \times_{X_0} X_1 \times_{X_0} \mathcal{L}_1'), \\
&\alpha^L: (\mathcal{L}_1 \times_{X_0} X_1 \times_{X_0} \mathcal{L}_1) \times_{\mathcal{L}_1 } (\mathcal{L}_1 \times_{X_0} X_1 \times_{X_0} \mathcal{L}_1') \rightarrow \mathcal{L}_1 \times_{X_0} X_1 \times_{X_0} \mathcal{L}_1',   \\
&\alpha^L((p,y,q) ,(q,y',r))= (p, yy',r), \\ 
&\forall ((p,y,q) ,(q,y',r))\in (\mathcal{L}_1 \times_{X_0} X_1 \times_{X_0} \mathcal{L}_1) \times_{\mathcal{L}_1 } (\mathcal{L}_1 \times_{X_0} X_1 \times_{X_0} \mathcal{L}_1'). \end{align*} 
 Moreover, by defining the 2-form $\Omega|_{(u,y,p)}=u^*\beta_1-p^*\beta_1'+{(x,\Phi'(p))}^*\omega_2-{(\Phi(u),y)}^*\omega_2$ on $\mathcal{L}_1 \times_{X_0} X_1 \times_{X_0} \mathcal{L}_1'$, we can see that this Morita equivalence is such that the quasi-symplectic forms are related as follows:
\begin{align*} &\left((\alpha^R)^* \Omega\right)|_{((u,y,v) ,(v,y',w))}= \Omega|_{(u,y,v)} + \Omega_{1,1}|_{(v,y',w)}- \omega_2|_{(x,x')}+ \omega_2|_{(y,y')}, \\ 
& x\Phi'(v)=\Phi(u)y,\, x'\Phi'(w)=\Phi'(v)y'  \\
&\forall ((u,y,v) ,(v,y',w))\in (\mathcal{L}_1 \times_{X_0} X_1 \times_{X_0} \mathcal{L}_1') \times_{\mathcal{L}_1' } (\mathcal{L}_1' \times_{X_0} X_1 \times_{X_0} \mathcal{L}_1'); \end{align*} 
and an analogous equation for $\alpha^L$.
\end{enumerate}  
\end{prop}
\begin{proof} In order to prove (1) notice that the groupoid structures over $\mathcal{L}_1,\mathcal{L}_1'$ are given by the maps $\alpha^R=\alpha^L$ for $\mathcal{L}_1'=\mathcal{L}$ and the fact that those maps define Morita equivalences of double Lie groupoids can be verified by an elementary inspection. In order to show that $\Omega_{1,1}$ is multiplicative over $X_1$ we can develop the equations $\delta \omega_2|_{((\Phi(u),\Phi(u'),z)}=0$, $\delta \omega_2|_{((x,\Phi(v),\Phi(v'))}=0$ and $\delta \omega_2|_{((\Phi(u),y,\Phi(v'))}=0$ and sum them to get that
\[ \omega_2|_{((\Phi(u),\Phi(u'))}-\omega_2|_{(\Phi(v),\Phi(v')}=\delta \widetilde{\Omega },   \] 
where $\widetilde{\Omega  }|_{(u,y,v)}=\omega_2|_{(x,\Phi(v))}-\omega_2|_{(\Phi(u),y)} $, and then the fact that $\delta \beta_1=\Phi^* \omega_2$ on $\mathcal{L}_1$ implies the equation $\delta \Omega_{1,1}=0$. Similar considerations imply $d \Omega_{1,1}=\delta \Omega_{1,0}$. Let's see that it indeed integrates the desired Dirac structure. By construction, $\mathcal{L}_1 \times_{X_0} X_1 \times_{X_0} \mathcal{L}_1 \rightrightarrows X_1$ is an action groupoid that integrates the action Lie algebroid $\mathtt{t}^* A_{L_1} \oplus \mathtt{s}^* A_{L_1} \rightarrow X_1$ associated to the map $\Gamma (A_{L_1}) \oplus  \Gamma (A_{L_1}) \rightarrow \mathfrak{X}(X_1) $ defined by $u\oplus v\mapsto \phi(u)^r+T\mathtt{i} (\phi(v))^l$ where $\phi$ is the Lie algebroid morphism induced by $\Phi$. So we just have to check the nondegeneracy of $\Omega_{1,1}$ and, by multiplicativity \cite[Lemma 3.1]{burint}, it is enough to verify it at $X_1$. Since $\Omega_{1,1}=\Omega $ for $\mathcal{L}_1=\mathcal{L}_1'$ we will check the nondegeneracy of the latter 2-form. So take $U\oplus V \in  A_{L_1}|_{\mathtt{t}(p) } \oplus  A_{L_1'}|_{\mathtt{s}(p)}$ and suppose that $\phi(U)^r|_p=-T \mathtt{i}( \phi'(V))^l|_p$, $\phi'$ being the Lie algebroid morphism determined by $\Phi'$. Take $W\in T_{\mathtt{t}( p)} X_1\cap \ker T_{\mathtt{t}( p)}\mathtt{s} \cap T_{\mathtt{t}( p)} \mathtt{t} $, then 
\begin{align*} &\langle i_{U\oplus V} \Omega_p,W \rangle =  \Omega((U,0,V),(0,W,0))=\\
&=\omega_2((\phi(U),0),(0  ,W^r))- \omega_2((\phi(U)^r|_p +T \mathtt{i}( \phi'(V))^l|_p,\phi'(V)),(W^r,0)=\\
&=\omega_2((\phi(U),0),(0 ,W^r))- \omega_2((0,\phi'(V)),(W^r,0)). \end{align*}  
Now let us observe that $\delta \omega_2|_{(p,\mathtt{s}(p),p^{-1} )}((0_p,\phi'(V),0_{p^{-1}}),(W^r,0_{\mathtt{s}(p) },0_{p^{-1}}))=0$ leads to 
\[ \omega_2(((0_p,\phi'(V)),(W^r,0_{\mathtt{s}(p) },0_{p^{-1}}))=-\omega_2((0_p,(\phi(U)^{-1})^l_{p^{-1}}),(W^r,0_{p^{-1}})) \]
and from $\delta \omega_2|_{(p,p^{-1},\mathtt{t}(p)  )}((0_p,0_{p^{-1}},\phi(U)^{-1}),(W^r,0_{\mathtt{s}(p) },0_{p^{-1}}))=0$ it follows that
\[ \omega_2((0_p,(\phi(U)^{-1})^l_{p^{-1}}),(W^r,0_{p^{-1}}))=\omega_2((0_{\mathtt{t}(p) },\phi(U)^{-1}),(W,0)). \] 
As a consequence, we get that
\[ \langle i_{U\oplus V} \Omega_p,W \rangle=\omega_2((\phi(U),0),(0 ,W^r))-\omega_2((0_{\mathtt{t}(p) },\phi(U)),(W,0)) \]
and this pairing gives us the induced metric at an isotropy Lie algebra which is nondegenerate. So if $i_{U\oplus V} \Omega_p|_{TpX_1} = 0$, then $\phi(U)=0$ and hence $\phi'(V)=0$. Since $X_1\rightrightarrows X_0$ is transitive, for every $Z\in T_{\mathtt{t}(p)}X_0$, there is $Z'\in A_{X_1}|_{\mathtt{t}(p)}$ which maps to it via the anchor. So the condition that $i_{U\oplus V} \Omega_p = 0$ gives us that
\[ \langle i_{U\oplus V} \Omega_p,(Z,(Z')^r_p,0)\rangle=\Omega_p((U,0,V),(Z,(Z')^r_p,0))=\beta_1(U,Z)=0. \]
But $\phi(U)=0$ and $U$ lies in the kernel of the map \eqref{eq:embcou} which is injective, thus $U=0$. Similarly, one shows that $V=0$. 

In order to conclude the proof that $\widetilde{D} (\Omega_{1,1} +\Omega_{1,0}+\Omega_{2,0})=0$ it remains to show that $\delta^V \Omega_{1,1}+\delta^H \Omega_{2,0}=0$, where we consider the groupoid structure over $\mathcal{L}_1$ as the vertical structure for consistency of notation. But notice that this equation is exactly the equation from item (2) if $\mathcal{L}_1=\mathcal{L}_1'$ and this last equation follows again from repeated application of $\delta \omega_2=0$. \end{proof}
\begin{rema} In the framework of shifted symplectic geometry, $\mathcal{L}_1 \times_{X_0} X_1 \times_{X_0} \mathcal{L}_1' \rightrightarrows X_1$ is the homotopic fibred product of the 2-shifted lagrangians $\mathcal{L}_1$ and $\mathcal{L}_1'$ in $X_1$ so it is a 1-shifted symplectic groupoid (quasi-symplectic groupoid) \cite{ptvv}, see also \cite{2lagpoi} for an interpretation and applications of this result in Poisson geometry. The Morita equivalence in this proposition is inspired by \cite{sevmorqua}. \end{rema} 
Now we will see how the previous result gives rise to a generalisation of the Lu-Weinstein construction of double symplectic groupoids associated to Lie bialgebroids \cite{luwei2}. Let $(A,B)$ be a Lie bialgebroid on $X_0$ such that the associated Courant algebroid $E=A\oplus B$ is transitive. Suppose that the Atiyah algebroid $A_{X_1}:=E/T^*X_0$ is integrable by a Lie groupoid $X_1 \rightrightarrows X_0$ and also $A$ and $B$ are respectively integrable by 2-shifted lagrangian groupoids $(\mathcal{L}(A) \rightrightarrows X_0,\Phi_A,\beta_1^A )$ and $(\mathcal{L}(B) \rightrightarrows X_0,\Phi_B,\beta_1^B)$ in $X_1 \rightrightarrows X_0$ with respect to some splitting of $\mathtt{t}^*E\oplus \mathtt{s}^* E$ defined by a closed 2-shifted 2-form $\omega=\omega_1+ \omega_2$. Then there is a double quasi-symplectic groupoid associated to $\mathcal{L}_1$ as in Proposition \ref{thm:2lagint}:
\[ \xymatrix{ \mathcal{G} :=\mathcal{L}(A) \times_{X_0} X_1 \times_{X_0} \mathcal{L}(A) \ar@<-.5ex>[r] \ar@<.5ex>[r]\ar@<-.5ex>[d] \ar@<.5ex>[d]& X_1  \ar@<-.5ex>[d] \ar@<.5ex>[d] & \mathcal{L}(B) \ar[l]_-{\Phi_B}\ar@<-.5ex>[d] \ar@<.5ex>[d] \\ \mathcal{L}(A) \ar@<-.5ex>[r] \ar@<.5ex>[r] & X_0 & X_0 \ar[l]_-{\text{id} }; }\] 
since $A$ and $B$ constitute a Lie bialgebroid with Courant algebroid $E$ and since $E$ covers $A_{X_1}$ surjectively, $\Phi_B$ is transverse to the anchor map of the multiplicative Dirac structure over $X_1 \rightrightarrows X_0$ corresponding to $\mathcal{G}$ and so we can pullback the double quasi-symplectic groupoid to a smooth  pullback (double) Lie groupoid $\Phi_B ^!\mathcal{G}$ over $\mathcal{L}(B)$. Let us observe that $\mathfrak{X}  =\Phi_B ^!\mathcal{G}$ can be described as follows:
\[ \mathfrak{X}=\{(x,v,u,y) \in \mathcal{L}(B) \times_{X_0} \mathcal{L}(A) \times \mathcal{L}(A) \times_{X_0} \mathcal{L}(B)| \, \Phi_B(x)\Phi_A(v)=\Phi_A(u)\Phi_B(y) \}; \]
this description makes manifest that $\mathfrak{X}$ is also a Lie groupoid over $\mathcal{L}(A)$. The source and target maps are the projections to $\mathcal{L}(A)$ (resp. $\mathcal{L}(B)$) and the two multiplications are a variation of \eqref{eq:douquamul} much like in the  original example \cite{luwei2}:
\begin{align*}  &\mathtt{m}^{\mathcal{L}(A)} ((x,v,u,y) ,(x',w,v,y') )=( xx',w,u, yy');\\
&\mathtt{m}^{\mathcal{L}(B)} ((x,v,u,y) ,(y,v',u',z) )=(x,vv',uu',z).
\end{align*} 
\begin{thm}\label{thm:dousym} The pullback groupoid $\mathfrak{X}$ is a double symplectic groupoid with sides $\mathcal{L}(A)$ and $\mathcal{L}(B) $ over $X_0$ integrating the Lie bialgebroid $(A,B)$ with
\begin{align*} \widehat{\Omega }_{(x,v,u,y)}= &(\Phi_B(x),\Phi_A(v))^*\omega_2- (\Phi_A(u),\Phi_B(y))^* \omega_2 + y^* \beta_1^B - x^* \beta_1^B + u^* \beta_1^A -v^* \beta_1^A \\
&\forall (x,v,u,y)\in \mathfrak{X}
 \end{align*}  
 as the symplectic form. \end{thm}   
\begin{rema} We can suppose that there is an additional Dirac structure $C \hookrightarrow A\oplus B$ and a 2-shifted lagrangian groupoid $\mathcal{L}(C) \rightrightarrows X_0$ that integrates it. So we can consider, more generally, the pullback groupoid $\Phi_B^!(\mathcal{L}(A) \times_{X_0} X_1 \times_{X_0} \mathcal{L}(C)) \rightrightarrows \mathcal{L}(B) $. This construction of quasi-symplectic groupoids is related to the problem of integrating the Poisson homogeneous spaces of $\mathcal{L}(B) \rightrightarrows X_0$ as a Poisson groupoid \cite{liuweixu2} and shall be discussed elsewhere. The insight that this construction of double symplectic groupoids could work in general goes back to \cite{sevlet}. \end{rema}
\begin{proof}[Proof of Theorem \ref{thm:dousym}] Notice that the pullback Dirac structure $L=\Phi_B^! \text{Lie} (\mathcal{G} )^H$ is a multiplicative Dirac structure inside the CA-groupoid $(\mathbb{T} \mathcal{L}(B) \rightrightarrows TX_0 \oplus B^*,\Phi_B^*(\omega_1 + \omega_2))$. But since $D \beta_1^B =\Phi_B^*(\omega_1 + \omega_2))$, applying the gauge transformation given by $-\beta_1^B$ to $L$ we get a Dirac structure inside the canonical CA-groupoid $\mathbb{T} \mathcal{L}(B) \rightrightarrows TX_0 \oplus B^*$ and so this Dirac structure is integrable by the corresponding gauge transformation of the pullback of $\Omega_{1,1}$ to $\mathfrak{X} $ and this gives us exactly $\widehat{\Omega } $. As a consequence, this 2-form is closed and multiplicative over $\mathcal{L}(B)$. Since we have that $\mathfrak{X}\cong\Phi_A^!(\mathcal{L}(B) \times_{X_0} X_1 \times_{X_0} \mathcal{L}(B))$, the same argument shows that $\widehat{\Omega } $ is multiplicative over $\mathcal{L}(A)$ (this also follows directly from the fact that $\mathcal{G} $ is a double quasi-symplectic groupoid). It remains to check the nondegeneracy of $\widehat{\Omega }$. We will see that $\widehat{\Omega } $ is nondegenerate over $X_0$ and, since it is multiplicative over both $\mathcal{L} (A)$ and $\mathcal{L}(B) $, it will have to be nondegenerate everywhere by \cite[Lemma 3.1]{burint}, indeed using this lemma once, we obtain that this 2-form is nondegenerate over the unit embedding of $\mathcal{L}(A)$. But then, using the fact that the 2-form is also multiplicative with respect to $\mathfrak{X}\rightrightarrows \mathcal{L}(A)$, we deduce the desired global nondegeneracy. To carry out this verification we shall use the fact that there is a canonical decomposition 
\[ T\mathfrak{X}|_{X_0}=TX_0\oplus C\oplus A\oplus B. \]
First, we will see that $\widehat{\Omega } $ induces an isomorphism between the core Lie algebroid $C$ and $T^*X_0$. Let $\phi_A$ and $\phi_B$ be the Lie algebroid morphisms induced by $\Phi_A$ and $\Phi_B$ respectively. Consider $c=(X,0,U,0)\in C_p$ and $W\in T_pX_0\subset T_p \mathfrak{X}  \hookrightarrow T_p(\mathcal{L}(B) \times \mathcal{L}(A) \times \mathcal{L}(A) \times \mathcal{L}(B))$, then
\[ \widehat{\Omega}_p(c,W)=\beta_1^A(U,W) - \beta_1^B(X,W)   \]
but notice that $\rank C=\dim X_0$ and the inclusions $A,B \hookrightarrow E $ are given by 
\[ a\mapsto \phi_A(a)\oplus -s_0^*\left(i_{(\phi_A(a),0)}\omega_2\right)+  \mathtt{u}^*(i_a \beta^A); \quad b\mapsto   \phi_B(b)\oplus -s_0^*\left(i_{(\phi_B(b),0)}\omega_2\right)+\mathtt{u}^*( i_b \beta^B )\quad    a\in A, b\in B; \]
with respect to the identification $E\cong \ker \mathtt{s}_{\omega }|_{X_0} $ associated with $\omega=\omega_1+ \omega_2$ as in \eqref{eq:embcou}. Since $A$ and $B$ have trivial intersection and $\phi_A(U)=\phi_B(X)$ 
we must have $\mathtt{u}^*(i_X \beta_1^B -i_U \beta_1^A)\neq 0$ if $c\neq 0$ and so our claim holds. On the other hand, $\widehat{\Omega } $ also defines an isomorphism $A\cong B^*$. In fact, take $u\in A_p$ and $v\in B_p$, then 
\begin{align}  &\widehat{\Omega }_p((\mathtt{a}_A(u),u,u,0 ),(v,0,\mathtt{a}_B(v),v )) = \omega_2((\mathtt{a}_A(u),\phi_A(u)),(\phi_B(v),0) )+ \label{eq:symfor} \\
&-\omega_2((\phi_A(u),0),(\mathtt{a}_B(v),\phi_B(v)) )+\beta_1^A(u,\mathtt{a}_B(v))- \beta_1^B(\mathtt{a}_A(u),v). \notag   
\end{align}
On the other hand, we have that the pairing of $E\cong \ker \mathtt{s}_{\omega }  $ restricted to $u,v$ is given by the canonical pairing of $\mathbb{T}X_1$ applied to $u$ and $v$ regarded as elements in $\ker \mathtt{s}_{\omega}$:
\begin{align}  &\langle \phi_A(u)\oplus -s_0^*\left(i_{(\phi_A(u),0)}\omega_2\right)+  \mathtt{u}^*(i_u \beta^A),\phi_B(v)\oplus -s_0^*\left(i_{(\phi_B(v),0)}\omega_2\right)+ \mathtt{u}^*( i_v \beta^B) \rangle; \label{eq:pai} 
\end{align} 
and substituting the decomposition $u=-u^{-1}+ \mathtt{a}_A(u)$ and $v=-v^{-1}+ \mathtt{a}_B(v)$ in \eqref{eq:pai} we get \eqref{eq:symfor} because $\mathtt{u}^*(i_u \beta^A)$ and $\mathtt{u}^*( i_v \beta^B)$ lie in $T^*X_0\hookrightarrow T^*X_1|_{X_0}$ identified with the annihilator of $\ker T \mathtt{t}$. Since $A$ and $B$ are in duality, $\widehat{\Omega } $ is nondegenerate when restricted to $A\oplus B$. Notice that the multiplicativity of $\widehat{\Omega}$
implies that it vanishes along $TX_0\oplus A$ and $TX_0\oplus B$. Thus the computation above together with the displayed decomposition of $T\mathfrak{X}|_{X_0}$ imply that $\widehat{\Omega}|_{X_0}$ is nondegenerate.

Since $\widehat{\Omega } $ determines a double symplectic groupoid structure, it differentiates to Poisson groupoid structures on $\mathcal{L}(A) \rightrightarrows X_0$ and $\mathcal{L}(B) \rightrightarrows X_0$. The LA-groupoid associated to, say $\mathcal{L}(A) \rightrightarrows X_0$, is $T^*\mathcal{L}(A) \rightrightarrows B$ by construction of $\mathfrak{X}$ and then these Poisson groupoid structures differentiate to the Lie bialgebroid $(A,B)$ as desired. \end{proof}
It appears that all the known examples of double symplectic groupoids are particular cases of the previous result. In most examples, the relevant transitive Courant algebroid is of the form $\mathbb{T}_H X_0 \oplus \mathfrak{g}$ where $\mathfrak{g}$ is a quadratic Lie algebra. The general way of showing that an integrable Dirac structure in a transitive Courant algebroid actually integrates to a 2-shifted lagrangian groupoid is by means of Van Est Theorem \cite{weivanest} as we shall discuss in \cite{2lagpoi}. However, in practice one can often find the desired integration by ad hoc considerations.
\begin{exa}[The Lu-Weinstein integration of a Lie bialgebroid over a point] The double symplectic groupoid that we get by applying Theorem \ref{thm:dousym} to a Lie bialgebroid over a point is the original construction of \cite{luwei2} since the 2-forms associated with a lagrangian subgroup of a Lie group equipped with a quadratic Lie algebra are zero. According to \cite{sevlet,sevmorqua}, this double symplectic groupoid is an example of a decorated moduli space of flat $G$-bundles. The double symplectic groupoids that appear in \cite{poigromod} are examples of this construction, see the example below. \end{exa} 
\begin{exa}[Integration of Lie bialgebroids in exact Courant algebroids] If the Courant algebroid  $E=A\oplus B$ of Theorem \ref{thm:dousym} is exact, then we can take $X_1 \rightrightarrows X_0$ as the pair groupoid over $X_0$ and $\Phi_A$ as the map $(\mathtt{t},\mathtt{s}):\mathcal{L}(A) \rightarrow X_0 \times X_0$ and the analogous definition for $\Phi_B$. In this case, the 2-shifted lagrangian groupoids $\mathcal{L}(A) \rightrightarrows X_0$ and $\mathcal{L}(B) \rightrightarrows X_0$ are quasi-symplectic groupoids themselves and the symplectic form $\widehat{\Omega  }$ on the double symplectic groupoid $\mathfrak{X}$ determined by Theorem \ref{thm:dousym} just depends on $\beta_1^A$ and $\beta_1^B  $ since $\omega_2$ vanishes. 

\begin{figure}
    \centering
    
   \begin{tikzpicture}[line cap=round,line join=round,>=stealth,x=1cm,y=1cm,scale=0.6]
\clip(0.8773230701015977,-0.451868545734714) rectangle (11.990433619994528,10.282385962684584);
\draw [line width=1pt] (6,5) circle (4.042190000482412cm);
\draw [line width=1pt] (6,5) circle (1.9544052803858265cm);
\draw [dashed,line width=1pt,<-] (7.42520238075843,3.662652560521199)-- (8.873503686657095,2.157065501498832);
\draw [dashed,line width=1pt,<-] (4.652621969355528,6.415723293068467)-- (3.1666147026899876,7.8829200053049435);
\draw [dashed,line width=1pt,->] (4.612936684202755,3.623135679171853)-- (3.1468033759400487,2.136685657415012);
\draw [dashed,line width=1pt,->] (7.261983888703337,6.492346027117439)-- (8.592647877584474,8.101205730495273);
\begin{scriptsize}
\draw[color=black] (8.2,3.2) node {$l$};
\draw[color=black] (3.7,7) node {$g$};
\draw[color=black] (3.7678543494929406,3.203390669223753) node {$h$};
\draw[color=black] (8.145746384299246,7) node {$k$};
\draw[color=black] (5.9,7.3) node {$b$};
\draw[color=black] (5.9,2.6) node {$b'$};
\draw[color=black] (3.7,5) node {$a$};
\draw[color=black] (8.3,5) node {$a'$};
\end{scriptsize}

\end{tikzpicture} 
\caption{The generators of the fundamental groupoid of a decorated surface with boundary.}
    \label{fig:decsur}
\end{figure}
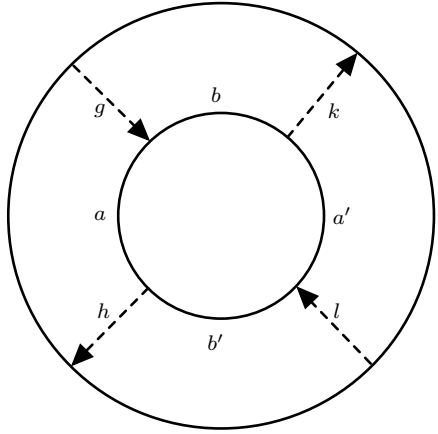

For instance, suppose that $(\mathfrak{g},\mathfrak{h},\mathfrak{k})$ is a Manin triple and $H \hookrightarrow G$, $K \hookrightarrow G$ are Lie subgroups corresponding to the Lie algebra inclusions
$\mathfrak{h}\hookrightarrow \mathfrak{g}$, $\mathfrak{k} \rightarrow \mathfrak{g}$, then we can construct action Lie groupoids $\mathcal{L}(A)=(K\times H )\times G$ and $\mathcal{L}(B)=(H\times K )\times G$, where it is assumed that the action in each case is the restriction of the action of $G\times G$ on $G$ defined by $(a,b)\cdot g=agb^{-1}$. The corresponding action Lie algebroids $A=(\mathfrak{k}\oplus \mathfrak{h})
\times G$ and $B=(\mathfrak{h}\oplus \mathfrak{k}) \times G$ are in duality by means of the metric on $\mathfrak{g}$ and the Courant algebroid that they induce is exact and it corresponds to the canonical Courant algebroid over $G$ twisted by the Cartan 3-form \cite{purspi}. In terms of the canonical splitting of this Courant algebroid \cite[\S 3.1]{purspi}, the 2-form that integrates the Dirac structure $A$ is
\[ \beta_1^A|_{(a,g,b)}=-\frac{1}{2}\langle h^{-1}dh\wedge (db)b^{-1}  \rangle +\frac{1}{2} \langle a^{-1}da \wedge (dg)g^{-1}  \rangle, \]
where $h=agb^{-1}$ and there is an analogous formula for $B$. 

There is a nice interpretation of the example discussed above as a moduli space of representations. Consider an oriented annulus $\Sigma$ with four marked points on each boundary component and let $V\subset \partial \Sigma$ be the set of all these marked points. Interestingly, the associated double symplectic groupoid $\mathfrak{X}$ can be identified with the subspace of groupoid morphisms in $\hom (\Pi_1(\Sigma,V),G)$ such that their values on the boundary arcs lie either in $A$ or in $B$ in an alternating fashion, 
Figure \ref{fig:decsur} illustrates this construction. The values of a representation $\rho \in \mathfrak{X}$ on the generators of $\Pi_1(\Sigma,V)$ labeled by $g,h,k,l$ lie in $G$ while the values of $\rho$ on the other indicated generators satisfy $a,a'\in K$ and $b,b'\in H$ and the constraints are $ka'l\in H$, $hb'l\in K$, $hag\in H$, $kbg\in K$. The symplectic form on $\mathfrak{X}$ coincides with the one determined by its structure as a decorated moduli space of representations as in \cite{sevmorqua}.

Other interesting examples of Lie bialgebroids in exact Courant algebroids come from generalised K\"ahler structures that we shall explore in future work. \end{exa}
\begin{exa}[Integration of Lie bialgebroids associated with Poisson actions] There is a natural family of Lie bialgebroids induced by Poisson actions of Poisson groups \cite[\S 4]{liuxubial}. Here we describe their double symplectic groupoids according to Theorem \ref{thm:dousym}. Suppose that $G$ is a complete connected Poisson group acting on an integrable Poisson manifold $(M,\pi)$ by a Poisson action. Let us assume that there is a dressing action of $G$ on $G^*$, the 1-connected integration of the dual Lie algebra $\mathfrak{g}^*$ and so the dressing actions endow $D:=G^* \times G$ with a Lie group structure that integrates the double Lie algebra $\mathfrak{d}=\mathfrak{g} \bowtie \mathfrak{g}^*$ \cite{luphd}. Consider the source-simply-connected integration $(\Sigma(M) \rightrightarrows M,\omega_\Sigma)$. There is a lifted Poisson action of $G$ on $\Sigma(M)$, denoted by $(a,x)\mapsto a\cdot x$ for all $a\in G$ and all $x\in \Sigma(M)$, and a moment map
\[ \mu:\Sigma(M) \rightarrow G^* \]
for this action which is a Lie groupoid morphism. Consider the lifted $G$-action on $\Sigma(M)$ together with the action of $\Sigma(M)$ on $G$ given by $a\cdot x=a^{\mu(x)}$ for all $a\in G$ and for all $x\in \Sigma(M)$, where $a^u$ denotes the dressing action of $u\in G^*$ on $a$ for all $u\in G^*$. Then we have the following property. The formula
\[ ((x,a),(y,b) )=(x(a\cdot y),(a^y)b) \]
for all suitable composable $x,y\in \Sigma(M)$ defines a Lie groupoid structure on the direct product $\Sigma(M) \times G$ that we denote by $\Sigma(M) \bowtie G \rightrightarrows M$; according to \cite[\S 9]{luliealg}, this Lie groupoid integrates the Dirac structure
\[ A=\{ \rho_p(u)+\pi^\sharp_p(\alpha)\oplus \alpha, u\oplus -\rho^*_p(\alpha):u\in \mathfrak{g},\alpha\in T^*_pM,\, p\in M \}\hookrightarrow \mathbb{T}M \oplus \mathfrak{d}. \] 
Note that $A$ was defined as just a Lie algebroid in \cite{luliealg}, its interpretation as a Dirac structure comes from \cite{quapoidir}. Let $\theta_K^l$ and $\theta_K^r$ be the left and right-invariant Maurer-Cartan forms on a Lie group $K$. We have that 
\[ \Psi: \mathcal{L}(A):=\Sigma(M) \bowtie G \rightarrow D, \quad \Psi(x,a)=(\mu(x),a) \]
is a Lie groupoid morphism and the 2-form
\[ \beta_1^A:= \Psi^* \langle \theta^l_{G^*}\wedge \theta^r_G \rangle- \text{pr}_1^*\omega_\Sigma  \]
endows the Lie groupoid morphism $\Phi_A=(\mathtt{t},\mathtt{s},\Psi):\mathcal{L}(A)\rightarrow X_1$ with a 2-shifted lagrangian structure, where $(X_1 \rightrightarrows X_0)= (M \times M \times D \rightrightarrows M)$ is the gauge groupoid associated to the trivial $D$-bundle over $M$ with the 2-shifted closed 2-form given by the pullback of \eqref{eq:bg} to $X_1 \rightrightarrows X_0$ (by considering $K=D$ in \eqref{eq:bg}). This follows from \cite[Thm. 4.8]{cuezhu} and shall also be explained in \cite{2lagpoi}. The Dirac structure $A$ admits a canonical complement inside $\mathbb{T}M \oplus \mathfrak{d}$:
\[ B=(TM \oplus 0)\oplus (0\oplus \mathfrak{g}^*). \]
So $A$ and $B$ constitute a Lie bialgebroid over $M$. We have that $B$ is integrable by $\mathcal{L}(B)$ given by the product of the pair groupoid over $M$ and $G^*$, $\mathcal{L}(B)$ is 2-shifted lagrangian in $X_1\rightrightarrows X_0$ with respect to the zero 2-form and $\Phi_B$ as the canonical inclusion. By applying Theorem \ref{thm:dousym} to this pair of 2-shifted lagrangian groupoids in duality we get the following double symplectic groupoid. Consider the product of symplectic manifolds $\Sigma(M) \times \overline{\Sigma(M)} \times D$, negating the sign of the second factor, and equip it with the product Lie groupoid structure over $\Sigma(M)\bowtie G$ of the pair groupoid over $\Sigma(M)$ and $D \rightrightarrows G$, let us call this the vertical groupoid structure. On the other hand, the groupoid structure of $\Sigma(M) \times \overline{\Sigma(M)} \times D$ over $\mathcal{L}(B)=M \times M \times G^* \rightrightarrows M$ is the one determined by the source and target maps
\begin{align*}
    &\mathtt{t}^H((x,y),(a,u))=(\mathtt{t}(x),\mathtt{t}(y),\mu(x)({}^au)\mu(y)^{-1}) \\  
    & \mathtt{s}^H((x,y),(a,u))=({a^{-1}}\cdot\mathtt{s}(x),{{(a^u)}^{-1}}\cdot\mathtt{t}(y),u) ;
\end{align*} 
for all $((x,y),(a,u))\in \Sigma(M) \times \Sigma(M) \times G \times G^* $, the horizontal multiplication map $\mathtt{m}^H$ for composable elements in $\Sigma(M) \times \overline{\Sigma(M)} \times D$ is given by
\[ (((x,y),(a,u)),((x',y'),(a',u')))\mapsto ((x(a\cdot x'),y({a^u}\cdot y')),((a^{\mu(x')})a',u')). \]

In particular, suppose that $\mathfrak{h}\hookrightarrow \mathfrak{g}$ is a Lie subalgebra such that $\mathfrak{h}^\circ\hookrightarrow \mathfrak{g}^*$ is an ideal and the quotient $\mathfrak{h}^*=\mathfrak{g}^*/\mathfrak{h}^\circ$ is an abelian Lie algebra. Then the action of $G^*$ on itself by left multiplication descends to a Poisson action of $G^*$ on $\mathfrak{h}^*=G^*/H^\circ$, where $H^\circ$ integrates $\mathfrak{h}^\circ$. Since $M=\mathfrak{h}^*$ inherits a linear Poisson structure, it is integrable by the cotangent groupoid $T^*H \rightrightarrows \mathfrak{h}^*$ and we can apply the previous construction to this example. The double symplectic groupoid that we get is the one constructed in \S5.2 of \cite{dynpoi}. 
\end{exa}     
\subsection{Double quasi-symplectic groupoids and local 2-shifted symplectic groupoids} Now we will see an interpretation of double quasi-symplectic groupoids in terms of shifted symplectic geometry with \cite{rajtan} as motivation. Since the following observation is not used in the rest of this paper, we shall be rather brief. For more details we refer to \cite{intlinf,kanrep,getdia,cuezhu}.

Let $(X_n)_{n\geq 0}$ be a simplicial manifold, we shall define its horn maps as
\begin{align*} &\lambda_{q,k}:X_q \rightarrow \Lambda_{q,k}=\{(x_0,\dots,x_{k-1},x_{k+1},\dots, x_q)\in X^q_{q-1}|d_ax_b=d_{b-1}x_a\, \forall a<b\} \\
& \lambda_{q,k}(x)=(d_0x,\dots,d_{k-1}x,d_{k+1}x,\dots,d_qx);
\end{align*} 
if $\lambda_{q,k}$ is a submersion for $0\leq k\leq q\leq l$, then $\Lambda_{l,k}$ is a manifold for all $k\leq l$ \cite[Lemma 2.4]{intlinf}.
\begin{defi}[\cite{kanrep}] A simplicial manifold $(X_n)_{n\geq 0}$ is a local Lie $n$-groupoid if $\lambda_{q,k}$ is a submersion for all $q\geq 1$ and an injective local diffeomorphism for all $q>n$. \end{defi} 
\begin{rema} If in the previous definition we require instead that $\lambda_{q,k}$ is a surjective submersion for all $q\geq 1$ and a diffeomorphism for all $q>n$ we get the notion of Lie $n$-groupoid \cite{intlinf}. These conditions are known as the {\em Kan conditions} while the ones above are the {\em local Kan conditions}. \end{rema}  
We shall see now how to transform a double Lie groupoid into a local Lie 2-groupoid. There is a functor that takes bisimplicial manifolds into simplicial manifolds which is called the codiagonal functor (also known as the bar construction \cite{rajtan} or classifying complex functor \cite{simp:goer}).
\begin{defi}[\cite{simp:goer}] Let $\mathfrak{X}= (X_{a,b})_{a,b\geq0}$ be a bisimplicial manifold. The {\em classifying functor} applied to $\mathfrak{X}$ is the simplicial manifold\footnote{It is assumed in this definition that the following fibred products exist. This is automatic in the case that interests us. } $\overline{W}\mathfrak{X}$:
\begin{align*} &(\overline{W}\mathfrak{X})_n=\{(x_n,x_{n-1},\dots,x_0)|x_a\in  X_{a,n-a},\, d_0^Hx_a=d^V_{a+1}x_{a+1}, \, \forall \, 0\leq a< n\} \\
& d_i(x_n,\dots,x_0)=(d_i^Vx_n,d_i^Vx_{n-1},\dots,d_i^Vx_{i+1},d_1^Hx_{i-1},d_2^Hx_{i-2},\dots,d_i^Hx_0) \\
& s_i(x_n,\dots,x_0)=(s_i^Vx_n,s_i^Vx_{n-1},\dots,s_i^Vx_i,s_0^Hx_i,s_{1}^Hx_{i-1},\dots,s_i^Hx_0). 
\end{align*}  \end{defi} 
If $\mathfrak{X}=( X_{ab})$ is the nerve of a double Lie groupoid, then $\overline{W}\mathfrak{X}  $ is a local Lie 2-groupoid \cite{rajtan}. In order to describe what is the result of applying $\overline{W} $ to a quasi-symplectic structure on $\overline{W}\mathfrak{X}  $ let us introduce the following concept. Let $X =(X_n)_{n\geq 0}$ be a Lie 2-groupoid. Its {\em tangent complex} $\mathcal{T}X$ is the 3-term complex of vector bundles over $X_0$ $(\mathcal{T}X,\partial) $:
\[ \left(\xymatrix{(\mathcal{T}X)_2  \ar[r]^-{\partial_2} &( \mathcal{T}X)_1  \ar[r]^-{\partial_1} & (\mathcal{T}X)_0 }\right)=\left(\xymatrix{\ker T\lambda_{2,2}|_{X_0} \ar[r]^-{Td_2} &\ker T\lambda_{1,1}|_{X_0} \ar[r]^-{-Td_1} & TX_0 }\right);\]
the convention is that $(\mathcal{T}  X)_i=\ker T\lambda_{i,i}|_{X_0}$ for all $i>0$ and $(\mathcal{T}  X)_i=0$ for negative $i$ but the only terms that are possibly nonzero are the ones above by the Kan condition. The cotangent complex of $X$ is the dual complex of $\mathcal{T}X$ so $(\mathcal{T}^*X)_i=(\mathcal{T}X)_{-i}^*$ for all $i$. The philosophy of shifted symplectic geometry is that an $n$-shifted symplectic form should be given by a closed $n$-shifted 2-form on $X$ that induces a quasi-isomorphism between $\mathcal{T}X $ and the $n$-shifted cotangent complex  $\mathcal{T}^*[n]X$ \cite{ptvv,getdia,bunger,cuezhu}. In particular, a normalised closed 2-shifted 2-form $\Omega =\Omega_2+\Omega_1+\Omega_0$ on $X$ induces a morphism of complexes
\begin{align} \lambda^{\Omega }:\mathcal{T}X \rightarrow \mathcal{T}^*[2]X \label{eq:getpai}    
\end{align} 
given by 
\begin{align*} &\lambda^{\Omega }_2(\xi)=-(s_1\circ s_0)^*i_\xi \Omega_2, \quad \lambda^{\Omega }_1(\eta)=s_0^* i_{Ts_1 (\eta)}\Omega_2-s_1^* i_{Ts_0 (\eta)}\Omega_2,  \quad \lambda^{\Omega }_0(\zeta)=i_{T (s_1\circ s_0)(\zeta)}{\Omega }_2, \\
& \forall \xi\in (\mathcal{T}X)_2,\, \forall\eta\in (\mathcal{T}X)_1, \, \forall \zeta\in (\mathcal{T}X)_0.  \end{align*}  
In this case, it is not difficult to see that equation $\delta \Omega_2=0$ implies that $\lambda^\Omega$ is indeed a morphism of complexes, for the general case see \cite{cuezhu}.
\begin{defi}[\cite{ptvv,getdia}]\label{def:2shisym} A normalised closed 2-shifted 2-form $\Omega $ on $X$ is symplectic if $\lambda^{\Omega }$ is a quasi-isomorphism. \end{defi} 
For instance, we have seen in Theorem \ref{thm:hetca} that the gauge groupoid associated to a heterotic Courant algebroid admits a 2-shifted symplectic groupoid structure.
\begin{thm}\label{thm:2sym} If $(\mathfrak{X},\Omega )$ is a double quasi-symplectic groupoid, then $(\overline{W} \mathfrak{X},\overline{W}^* \Omega )  $ is a 2-shifted symplectic local Lie 2-groupoid which is a Lie 2-groupoid if $\mathfrak{X} $ is a full double Lie groupoid. \end{thm}

The proof of this theorem is based on two observations: (1) $\overline{W}$ induces a morphism of total complexes of differential forms and (2) the horizontal nondegeneracy of a closed 2-form implies the nondegeneracy of the 2-form that results from applying $\overline{W} $.
\begin{lem}\label{lem:biscom} If $\mathfrak{X}=(X_{b,c})$ is a bisimplicial manifold, then the natural pullback map 
\[ \overline{W}^*_n:\bigoplus_{a+b+c=n}\Omega^a(X_{b,c}) \rightarrow \bigoplus_{p+q=n} \Omega^p((\overline{W} \mathfrak{X})_q )   \] 
satisfies that $\overline{W}^*_{n+1} \widetilde{D} = D \overline{W}^*_n$ for all $n$. \end{lem}
\begin{proof} We just have to check that our claim holds on $\Omega^i(X_{j,k})$. The projections $\text{pr}_j:(\overline{W} \mathfrak{X})_{q=b+c}\rightarrow  X_{j,q-j}$ satisfy that
\begin{align*} \text{pr}_j\circ d_i= \begin{cases} d_i^V\circ \text{pr}_{j+1}, \quad \text{if $j>i-1$}   \\
d_{i-j}^H\circ \text{pr}_j, \quad \text{if $0\leq j \leq i-1$}.  \end{cases}    \end{align*}   
As a consequence, we have that
\[ \delta \overline{W}^*= \sum_{i=0}^n (-1)^i d_i^* \text{pr}_j^*=\sum_{i<j+1} (-1)^i\text{pr}_{j+1}^* (d_i^V)^*+ \sum_{i\geq j+1} (-1)^i\text{pr}_j^* (d_{i-j}^H)^*. \]
But we have that $d_0^H \circ \text{pr}_{j} =d_{j+1}^V\circ \text{pr}_{j+1}$, therefore
\[ \delta \overline{W}^*=\text{pr}_{j+1}^*\sum_{ i=0}^{ j+1} (-1)^i (d_i^V)^*+ (-1)^j\text{pr}_j^*\sum_{ i=0}^{n-j} (-1)^i (d_{i}^H)^*= \text{pr}_{j+1}^* \delta^V +(-1)^j\text{pr}_{j}^*\delta^H .  \] 
On the other hand, the desired relation clearly holds for the de Rham component of the differentials so the result holds. \end{proof} 
\begin{lem} If $(\mathfrak{X},\Omega )$ is a double quasi-symplectic groupoid, then $\lambda^{\overline{W}^* \Omega }$ establishes a quasi-isomorphism between $\mathcal{T} \overline{W}  \mathfrak{X}$ and $\mathcal{T}^*[2] \overline{W}  \mathfrak{X}$. \end{lem} 
\begin{proof} As we shall see, the Kan conditions imply that $ \mathcal{T} \overline{W}  \mathfrak{X}$ has only three nonzero terms:
\begin{align*}  (\mathcal{T} \overline{W}  \mathfrak{X})_2=\left\{(0,v,w)\in T (\overline{W}\mathfrak{X})_2|_{X_{0,0}}  \left| \begin{aligned}  &T\mathtt{s}^V(v)=0,\, 
T\mathtt{s}^H(v)=0,\\ 
&w_1=T\mathtt{i}T\mathtt{t}^V(v),\,\, w=(w_1,w_2)   \end{aligned} \right. \right\}\subset T (X_{2,0} \times X_{1,1} \times X_{0,2})\end{align*}
is identified with the core $C$ of the horizontal LA-groupoid induced by $(\mathfrak{X},\Omega )$. The differential is $\partial_2 (0,v,w)=(T d_1^H (v),T d_2^H(w))=(\mathtt{a}^H(v),T\mathtt{i}(\mathtt{a}^V(v)))$, where $\mathtt{a}^H,\mathtt{a}^V$ are the respective anchors of the horizontal and vertical LA-groupoids. Next we have 
\[ (\mathcal{T} \overline{W}  \mathfrak{X})_1=\left\{ (x,y)\in  T (\overline{W}\mathfrak{X})_1|_{X_{0,0}}| T d_0^V(x)=0  \right\}\subset T X_{1,0} \times TX_{0,1} \]
that we can identify with the vector bundle $A_{1,0}\oplus A_{0,1}$ and the differential $\partial_1:(\mathcal{T} \overline{W}  \mathfrak{X})_1 \rightarrow TX_{0,0}$ is the sum of the anchors with a negative sign. With respect to this identification, $\partial_2$ becomes $v\mapsto \mathtt{a}^H(v)\oplus -\mathtt{a}^V(v)$ for all $v\in C$. 

The LA-groupoid $A^H \rightrightarrows A_{0,1}$ sits inside $(\mathbb{T}X_{1,0} \rightrightarrows TX_{0,0}\oplus A_{1,0}^*,\Omega_{2,0}+ \Omega_{1,1})$ as a multiplicative Dirac structure (Prop. \ref{pro:difdir}), in particular, the map $u\mapsto (T \mathtt{t}^H(u),i_u \Omega_{1,1})$ is an inclusion over $X_{0,0}$. It turns out that this fact together with the equality $\rank A^H=\dim X_{1,0}$ imply the nondegeneracy of $\overline{W}^* \Omega $. In order to see that, let us first describe the components of $\lambda^{\overline{W}^* \Omega}$:
\begin{align*} &\langle \lambda^{\overline{W}^* \Omega}_{2}(0,v,w), a \rangle =- \Omega_{1,1}(v, T (s_0^H\circ s_0^V) (a) ), \quad \forall (0,v,w)\in (\mathcal{T} \overline{W}  \mathfrak{X})_2|_p,\, a\in T_pX_{0,0},   \\ 
& \langle \lambda^{\overline{W}^* \Omega}_1|_p(x,y),(x',y') \rangle=\langle (s_0^V)^* i_{Ts_1^V (x)}\Omega_{2,0},x' \rangle - \langle (s_1^V)^* i_{Ts_0^V (x)}\Omega_{2,0}, x' \rangle + \\
& +\Omega_{1,1}(T s_0^H(x), T s_0^V(y'))-\Omega_{1,1}(T s_0^V(y), T s_0^H(x')),\quad \forall (x,y),(x',y')\in (\mathcal{T} \overline{W}  \mathfrak{X})_1|_p  \\  
& \lambda^{\overline{W}^* \Omega}_0=-(\lambda^{\overline{W}^* \Omega}_2)^*.  
  \end{align*} 
On the other hand, we know that the maps 
\begin{align*} &P:C \rightarrow A_{1,0}\oplus T^*X_{0,0}, \quad v\mapsto (Td_1^H(v),-(s_0^H\circ s_0^V)^*i_v \Omega_{1,1}) \\
&Q:A_{0,1} \rightarrow TX_{0,0}\oplus A_{1,0}^*, \quad y\mapsto (-Td_1^H(y),-(s_0^H)^* i_{y} \Omega_{1,1} )  \end{align*}  
are injective since they correspond to the inclusion of $A^H|_{X_{0,0}}$ as a Dirac structure (up to sign).

Since $\widetilde{D} \Omega =0$, Lemma \ref{lem:biscom} implies that $D\overline{W}^* \Omega =0 $ and so we have a morphism of complexes $\lambda^{\overline{W}^* \Omega}:\mathcal{T} \overline{W}  \mathfrak{X} \rightarrow \mathcal{T}^*[2] \overline{W}  \mathfrak{X}$. This morphism is a quasi-isomorphism if and only if the induced mapping cone is exact and we can arrange it as follows:
\begin{align*}  &\xymatrix{0 \ar[r] & C \ar[r]^-{F } &  A_{1,0}\oplus T^*X_{0,0}\oplus A_{0,1} \ar[r]^-{G} & TX_{0,0}\oplus A_{1,0}^*\oplus A_{0,1}^*  \ar[r]^-{H} & C^* \ar[r] &0 } \\
&F=(P,-\mathtt{a}^V) ,\quad G=\begin{pmatrix} Q' & \partial_1^*| & Q \end{pmatrix}   \quad H= \begin{pmatrix} \lambda^{\overline{W}^* \Omega}_0 & -(\mathtt{a}^H)^* & (\mathtt{a}^V)^* \end{pmatrix}; \end{align*}
where $Q'$ is an expression such that $\text{pr}_{A_{0,1}^*}\circ (Q',\partial_1^*)=-Q^*$. Notice that this sequence is self-dual and $F$ is injective since it has $P$ as one of its components, hence $ H=-F^*$ is surjective. Let us denote $\rank C=m$, $\rank A_{0,1}=n$, $\rank A_{1,0}=r$ and $\rank TX_{0,0}=s$. Since $m+n=r+s$, we have that $\rank (\ker H)=2n $ but by analysing the projections we get that $\rank G=\rank Q+ \rank \text{pr}_{A_{0,1}^*}\circ (Q',\partial_1^*)=2n$ and so the sequence is exact. \end{proof}
\begin{proof}[Proof of Theorem \ref{thm:2sym}] It only remains to show that the Kan conditions hold in this situation from which it follows that the associated tangent complex has length three. But this follows from \cite[Thm. 4.5]{rajtan}. \end{proof}  
\appendix
\section{On the structure maps of an exact twisted CA-groupoid}
In this appendix we will continue using the notation and assumptions of Proposition \ref{thm:excagpd} and its proof.
\begin{prop}\label{pro:twimul}
    The Dirac structure $L_\omega$ of \eqref{eq:twimul} is the graph of a groupoid multiplication on $\mathbb{E}_1\rightrightarrows \mathbb{E}_0$.
\end{prop}
\begin{proof} {\em Definition of the source and target maps.} The deformed source and target maps $\mathtt{s}_{\omega },\mathtt{t}_{\omega }:\mathbb{E}_1 \rightarrow \mathbb{E}_0$ are forced to be defined by the equations:
\begin{align*} &\mathtt{m}_{G_1}^* \alpha =\text{pr}_1^* \alpha + \text{pr}_2^*(p_{T^*}(  \mathtt{s}_{\omega }(X\oplus \alpha)) )-i_{(X,T\mathtt{s}_{G_1}(X) ) }\omega_2  \\ 
& \mathtt{m}_{G_1}^* \alpha =\text{pr}_1^*(p_{T^*}(  \mathtt{t}_{\omega }(X\oplus \alpha)) )+\text{pr}_2^* \alpha  -i_{(T\mathtt{t}_{G_1}(X),X ) }\omega_2;\notag \end{align*} 
for all $X\oplus \alpha \in \mathbb{T}_gG_1 $, where $\text{pr}_i:G_2 \rightarrow G_1$ are the projections. We can see that that the normalisation condition for $\omega_2 $ implies that these equations induce well defined components in $A^*$ for the source and target maps while $p_T \circ \mathtt{s}_{\omega }=T \mathtt{s}_{G_1}$ and $p_T \circ \mathtt{t}_{\omega }=T \mathtt{t}_{G_1}$ are forced because the anchor $p_T$ has to  be a VB-groupoid morphism. Therefore, in terms of $\mathtt{s}_{T^*G}$ and $ \mathtt{t}_{T^*G} $, the canonical source and target maps of $T^*G_1$, we have that:
\begin{align} &\langle \mathtt{s}_{\omega } (X\oplus \alpha ), u \rangle= \langle \mathtt{s}_{T^*G} (\alpha)  ,u \rangle +\langle i_{(X,T \mathtt{s}(X) )} \omega_2,(0_g,u) \rangle    \label{eq:prstmap} \\
& \langle \mathtt{t}_{\omega } (X\oplus \alpha ), w \rangle= \langle \mathtt{t}_{T^*G}( \alpha) ,w \rangle +\langle i_{(T \mathtt{t}(X),X )} \omega_2,(w,0_g) \rangle \notag 
\end{align} 
for all $X\oplus \alpha \in \mathbb{T}_gG_1 $ and all $u\in A_{\mathtt{s}(g) }\cong \ker T_{\mathtt{s}(g) } \mathtt{t}  $, $w\in A_{\mathtt{t}(g) }\cong \ker T_{\mathtt{t}(g) } \mathtt{s}$. 

{\em Definition of the multiplication}. The multiplication $\mathtt{m}_{\omega  }:\mathbb{E}_2= \mathbb{E}_1 \times_{\mathbb{E}_0} \mathbb{E}_1 \rightarrow \mathbb{E}_1  $ is defined by its graph $L_{\omega } \subset \overline{\mathbb{T}G_1} \times  \mathbb{T}G_1 \times \mathbb{T}G_1$. In fact, equation $\mathtt{m}^*_{G_1} \alpha =\text{pr}_1^*\beta + \text{pr}_2^*\gamma-i_{(Y,Z)}\omega_2$ defines $\alpha =p_{T^*}(\mathtt{m}_{\omega }(Y\oplus \beta ,Z\oplus \gamma))$ uniquely for every $Y\oplus \beta $ and $Z\oplus \gamma $ composable in $\mathbb{E}_1 $. Since $\mathtt{m}_{G_1}$ is a submersion, we just have to check that $\theta:=\text{pr}_1^*\beta + \text{pr}_2^*\gamma-i_{(Y,Z)}\omega_2$ is $\mathtt{m}_{G_1}$-basic, meaning that $\theta(u,v)=0$ for all $(u,v)\in \ker T \mathtt{m}_{G_1}$. If $T \mathtt{m}_{G_1}(u,v)=0$, we can write $u=T \mathtt{m}_{G_1}(0,u^l)$ and $v=T \mathtt{m}_{G_1}(v^r,0)$, where $u^l=-U\oplus \mathtt{a}(U)\in A\oplus TG_0$ and $v^r=U\oplus 0\in A\oplus TG_0$, since we must also have $T \mathtt{m}_{G_1}(u^l,v^r)=0$. Using \eqref{eq:prstmap} we have that 
\begin{align*}  &\langle \beta ,u \rangle =\langle \mathtt{m}_{G_1}^* \beta ,(0,u^l) \rangle =\langle \text{pr}_1^* \beta + \text{pr}_2^*(p_{T^*}(  \mathtt{s}_{\omega }(Y\oplus \beta )) )-i_{(Y,T\mathtt{s}_{G_1}(Y) ) }\omega_2 , (0,u^l)\rangle \\ 
& \langle \gamma,v \rangle = \langle \mathtt{m}_{G_1}^* \gamma , (v^r,0) \rangle  =\langle \text{pr}_1^*(p_{T^*}(  \mathtt{t}_{\omega }(Z\oplus \gamma )) )+\text{pr}_2^* \gamma  -i_{(T\mathtt{t}_{G_1}(Z),Z ) }\omega_2,(v^r,0) \rangle ; \end{align*}     
and the composability of $Y\oplus \beta $ and $Z\oplus \gamma $ means that $ p_{T^*}(  \mathtt{s}_{\omega }(Y\oplus \beta  ))= p_{T^*}(  \mathtt{t}_{\omega }(Z\oplus \gamma )) $ and since these covectors live in the conormal bundle of $G_0$ we have that 
\[ \langle p_{T^*}(  \mathtt{s}_{\omega }(Y\oplus \beta  )), u^l \rangle =- \langle (p_{T^*}(  \mathtt{t}_{\omega }(Z\oplus \gamma )), v^r \rangle \] 
by decomposing $u^l$ and $v^r$ as above. But then we have that  
\begin{align}  \langle  \theta  ,(u,v) \rangle =-\langle i_{(Y,Z)}\omega_2,(u,v)\rangle -\langle i_{(T\mathtt{t}_{G_1}(Z),Z ) }\omega_2,(v^r,0)\rangle-\langle i_{(Y,T\mathtt{s}_{G_1}(Y) ) }\omega_2 , (0,u^l)\rangle. \label{eq:mulmppr} \end{align} 
On the other hand, using again the decomposition of $u^l$ and $v^r$ we have that  
\[ -\omega_2 (Y,T\mathtt{t}_{G_1}(Z)),(u,v^r))=-\omega_2 (Y,T\mathtt{t}_{G_1}(Z)),(u,T\mathtt{t}_{G_1}(u) ))+\omega_2 (Y,T\mathtt{t}_{G_1}(Z)),(0,u^l)); \] 
but $\omega_2 (Y,T\mathtt{t}_{G_1}(Z)),(u,T\mathtt{t}_{G_1}(u) ))=0$ by the normalisation condition. The facts that $T \mathtt{m}_{G_1} (u,v^r)=0$, $\delta \omega_2=0$ and
\begin{align*} &0=\delta \omega_2((Y,T\mathtt{t}_{G_1}(Z),Z) ,(u,v^r,0))=\omega_2 (Y,T\mathtt{t}_{G_1}(Z)),(u,v^r))-\omega_2((Y,Z),(u,v))+ \\ 
& +\omega_2((Y,Z ),(T \mathtt{m}_{G_1} (u,v^r),0))-\omega_2((T\mathtt{t}_{G_1}(Z),Z),(v^r,0)), \end{align*} 
imply that, by substituting in \eqref{eq:mulmppr}, $\langle \theta ,(u,v) \rangle =0$ as desired. As a consequence, $L_{\omega }$ is really the graph of a vector bundle map $\mathtt{m}_{\omega }: \mathbb{E}_2 \rightarrow \mathbb{E}_1$ that covers $T \mathtt{m}_{G_1}$ as desired. By calculations similar to those perfomed above, we get that $\mathtt{t}_{\omega }\circ \mathtt{m}_{\omega }=\mathtt{t}_{\omega }\circ \text{pr}_1$ and $\mathtt{s}_{\omega }\circ \mathtt{m}_{\omega }=\mathtt{s}_{\omega }\circ \text{pr}_2$. By the definition of the source and target maps, $\mathtt{m}_{\omega }$ satisfies the unit axiom. On the other hand, equations $\mathtt{s}_{\omega }\circ \mathtt{u}_{\omega }= \text{id}_{\mathbb{E}_0 }= \mathtt{t}_{\omega }\circ \mathtt{u}_{\omega }$ follow directly from the definition \eqref{eq:prstmap} and the normalisation of $\omega_2$.

{\em Claim: we have that $\delta \omega_2=0$ is equivalent to the associativity for $\mathtt{m}_{\omega } $}. This can be checked as follows. Let us consider the four different face maps $d_i:\mathbb{E}_3 \rightarrow \mathbb{E}_2$ composed with the inclusion $\mathbb{E}_2 \cong L_{\omega }\subset \overline{\mathbb{E}}_1 \times \mathbb{E}_1 \times \mathbb{E}_1   $, $F_i: \mathbb{E}_3 \rightarrow \overline{\mathbb{E}}_1 \times \mathbb{E}_1 \times \mathbb{E}_1$:
\begin{align*}  F_0(\xi_1,\xi_2,\xi_3)=(\mathtt{m}_{\omega }(\xi_1,\xi_2),\xi_1,\xi_2),\quad F_1(\xi_1,\xi_2,\xi_3)=(\mathtt{m}_{\omega }(\mathtt{m}_{\omega } (\xi_1,\xi_2),\xi_3),\mathtt{m}_{\omega } (\xi_1,\xi_2),\xi_3) \\
F_2(\xi_1,\xi_2,\xi_3)=(\mathtt{m}_{\omega }(\xi_1,\mathtt{m}_{\omega } (\xi_2,\xi_3)),\xi_1,\mathtt{m}_{\omega } (\xi_2,\xi_3)),\quad  F_3(\xi_1,\xi_2,\xi_3)=(\mathtt{m}_{\omega }(\xi_2,\xi_3),\xi_2,\xi_3);
\end{align*} 
we can see the graph of the composition of $F_1$ with $F_0 \times \text{id}$ as a Courant morphism over the graph of the map $G_3 \rightarrow G_1$ defined by $(a,b,c)\mapsto \mathtt{m}_{G_1} (\mathtt{m}_{G_1}(a,b),c)$; analogously $F_2\circ (\text{id} \times F_3)$ is a Courant morphism over the map $(a,b,c)\mapsto \mathtt{m}_{G_1}(a,\mathtt{m}_{G_1}(b,c)) $ from $G_3 $ to $G_1$. The projections 
\[ p_{T^*}(\mathtt{m}_{\omega }(\mathtt{m}_{\omega } (\xi_1,\xi_2),\xi_3)), \quad p_{T^*}(\mathtt{m}_{\omega }(\xi_1,\mathtt{m}_{\omega } (\xi_2,\xi_3))) \] 
are equal if and only if $F_1\circ (F_0 \times \text{id})= F_2\circ (\text{id} \times F_3)$ and, by spelling out the defining equations, this is equivalent to the fact that
\begin{align*} & \omega_2((X_1,X_2),(Y_1,Y_2))+\omega_2((T\mathtt{m}_{G_1}(X_1,X_2),X_3),(T\mathtt{m}_{G_1}(Y_1,Y_2),Y_3))=\\
& \omega_2((X_1,T\mathtt{m}_{G_1}(X_2,X_3)),(Y_1,T\mathtt{m}_{G_1}(Y_2,Y_3))+\omega_2((X_2,X_3),(Y_2,Y_3)), \end{align*}   
where $X_i=p_T(\xi_i)$ and $(Y_1,Y_2,Y_3)\in TG_3$.

{\em Definition of the inversion $\mathtt{i}_{\omega }:\mathbb{E}_1 \rightarrow \mathbb{E}_1$}. We have to see that, for every $X\oplus \alpha \in \mathbb{E}_1$, the equation
\[ \mathtt{m}_{G_1}^*( p_{T^*} (\mathtt{t}_{\omega }(X\oplus \alpha ) ))=\text{pr}_1^* \alpha + \text{pr}_2^*(p_{T^*}(\mathtt{i}_{\omega } (X\oplus \alpha )))-i_{(X,T\mathtt{i}_{G_1}(X))}\omega_2  \]
has a (unique) solution for $p_{T^*}(\mathtt{i}_{\omega } (X\oplus \alpha ))$; in other words, we have to verify that 
\[ \eta:=\mathtt{m}_{G_1}^*( p_{T^*} (\mathtt{t}_{\omega }(X\oplus \alpha ) ))-\text{pr}_1^* \alpha +i_{(X,T\mathtt{i}_{G_1}(X))}\omega_2 \]
vanishes on vectors of the form $(Y,0)\in TG_2$. But we have that 
\[ \langle \mathtt{m}_{G_1}^*( p_{T^*} (\mathtt{t}_{\omega }(X\oplus \alpha ) )),(Y,0) \rangle =\langle \mathtt{m}_{G_1}^* \alpha -\text{pr}_2^*\alpha+i_{(T \mathtt{t}_{G_1}(X),X )}\omega_2, (Y^r,0) \rangle; \]
where $Y^r=T \mathtt{m}_{G_1}(Y,0)$. It follows that 
\[ \langle \eta,(Y,0) \rangle =\langle i_{(X,T\mathtt{i}_{G_1}(X))}\omega_2,(Y,0) \rangle +\langle i_{(T \mathtt{t}_{G_1}(X),X )}\omega_2, (Y^r,0) \rangle \]
and by developping the equation $\delta \omega_2((X,T\mathtt{i}_{G_1}(X),X),(Y,0,0)=0 $ we get that indeed $\langle \eta,(Y,0) \rangle=0$. Similarly, $\mathtt{t}_{\omega }\circ \mathtt{i}_{\omega }=\mathtt{s}_{\omega }$ and  $\mathtt{s}_{\omega }\circ \mathtt{i}_{\omega }=\mathtt{t}_{\omega }$. \end{proof}
\subsection*{Acknowledgements} The author thanks H. Bursztyn, M. Cueca and E. Meinrenken for enlightening conversations related to this work. The author is also grateful to the referees for their valuable comments and suggestions.
\printbibliography
\end{document}